\documentclass[BCOR8mm,DIV14]{scrartcl}
\usepackage{umlaut}     
\usepackage{bbm}  
\usepackage[english]{babel}
\usepackage{amsmath,amssymb,amsthm}
\newtheorem{theorem}{Theorem}  

\newtheorem{corollary}{Corollary}
\newtheorem{proposition}{Proposition}
\newtheorem{lemma}{Lemma}
\newtheorem{remark}{Remark}
\newtheorem{definition}{Definition} 
\begin{document}
\title{ Angles and Polar Coordinates In Real Normed Spaces  }
\author{VOLKER  THÜREY  \\  Rheinstr. 91  \\   28199 Bremen,  Germany     
     \thanks{T: 49 (0)421/591777 \  \   E-Mail: volker@thuerey.de }   }
 \maketitle
  \centerline{ 	  MSC-class: 52A10   \quad  \quad  
            Keywords: angles, normed space, polar coordinates    }   
  \begin{abstract}
     \centerline{Abstract}     
     We try to  create a  wise  definition of 'angle spaces'. 
     Based on an idea of Ivan Singer, we  introduce a new concept of an angle in real Banach spaces, 
     which  generalizes the  euclidean angle in  Hilbert spaces.  
     With this angle it is shown that in  every two-dimensional subspace of a real Banach space 
     we can describe elements uniquely by polar coordinates. 
     \end{abstract} 
	 \tableofcontents 
	  	    \section{Introduction}      
	  In a real  inner product space  \ $ ( X , < \ .. \ | \ ..> ) $ \ it is well-known that 
	   the inner product can be expressed by the norm, namely for \ \  $ \vec{x} , \vec{y} \in X $ , \ \  
	    $ \vec{x} \neq \vec{0} \neq \vec{y} $ , 
	    \begin{equation*}
         <\vec{x} \ | \ \vec{y}> \  = \ 
       \frac{1}{4} \cdot ( \: \|\vec{x}+\vec{y}\|^{2} -   \|\vec{x}-\vec{y}\|^{2} \: ) \  =  \    
          \frac{1}{4} \cdot  \|\vec{x}\| \cdot \|\vec{y}\| \cdot 
         \left[  \ \left\| \frac{\vec{x}}{\|\vec{x}\|}   +  \frac{\vec{y}}{\|\vec{y}\|} \right\|^{2} \  - \
         \left\| \frac{\vec{x}}{\|\vec{x}\|}   -   \frac{\vec{y}}{\|\vec{y}\|} \right\|^{2} \  \right] \ .
      \end{equation*}
    Furthermore we have  for all \  $ \vec{x},\vec{y} \neq \vec{0} $ \ the  euclidean angle  
     \begin{equation*} 
   \angle_{Euclid} (\vec{x}, \vec{y}) := \arccos{\frac{< \vec{x} \: | \: \vec{y} >}{ \|\vec{x}\| \cdot \|\vec{y}\|}} 
   \ = \  \arccos{    \left( \ \frac{1}{4} \cdot   
         \left[  \ \left\| \frac{\vec{x}}{\|\vec{x}\|}   +  \frac{\vec{y}}{\|\vec{y}\|} \right\|^{2} \  - \
      \left\| \frac{\vec{x}}{\|\vec{x}\|}   -   \frac{\vec{y}}{\|\vec{y}\|} \right\|^{2} \  \right] \  \right) } \ , 
      \end{equation*} 
   which is defined in terms of the norm, too.   \\  
    In this  paper we deal with real topological vector spaces \ $X$ \ provided with a continuous  map \ 
     $\|..\|  \longrightarrow \  \mathbbm{R}^{+}  \cup  \{ 0  \}$ \ which is absolute homogeneous, i.e. \ 
     \ $ \| r \cdot \vec{x} \| = |r| \cdot \|\vec{x}\| $ \  for \  $ \vec{x} \in X \ , \ r \in   \mathbbm{R} $. \  
    We call such pairs \ $ ( X, \|..\| ) $ \ homogeneously weighted vector spaces  (or {\bf hw spaces}).  
    The subset \ ${ \sf Z} \ :=  \ \{ \vec{x} \in X  \ | \ \|\vec{x}\| = 0 \ \} \ \subset X $ \ 
       is  called the {  zero-set} \ of ( $ X , \|..\| $ ) .  \\
    Following the lines of an inner product we define for such spaces a product   
   \quad   $ <.. \: | \: .. >_\spadesuit \ : \  X^{2} \longrightarrow   \mathbbm{R} $ ,    \ \ 
      writing  for all \  $  \vec{x}, \vec{y} \in X  $ : \qquad \
    \begin{equation*}   
         < \vec{x} \: | \:  \vec{y} >_\spadesuit   \ \  :=  \ \             
             \begin{cases}
           0  \quad  &
       \quad \mbox{for} \quad  \|\vec{x}\| \cdot \|\vec{y}\|  =     0   \\    
        \frac{1}{4} \cdot  \|\vec{x}\| \cdot \|\vec{y}\| \cdot 
         \left[  \ \left\| \frac{\vec{x}}{\|\vec{x}\|}   +  \frac{\vec{y}}{\|\vec{y}\|} \right\|^{2} \  - \
         \left\| \frac{\vec{x}}{\|\vec{x}\|}   -   \frac{\vec{y}}{\|\vec{y}\|} \right\|^{2} \  \right]
            &  \quad \mbox{for} \quad   \|\vec{x}\| \cdot \|\vec{y}\|  \neq    0  \ \ \ . \\            
                         \end{cases} 
     \end{equation*} 
    and it is easy to show that such product fulfils  the symmetry (  $ <  \vec{x} | \vec{y} >_\spadesuit  \
    = \ <  \vec{y} | \vec{x} >_\spadesuit  $ ), the positive semidefiniteness 
   ( $  <  \vec{x} | \vec{x} >_\spadesuit  \ \geq \ 0 $ ), and  the homogenity 
   ( $  < r \cdot \vec{x} | \vec{y} >_\spadesuit  \ = \ r \cdot <  \vec{x} | \vec{y} >_\spadesuit $), 
   for \ $  \vec{x},  \vec{y} \in X $ , \  $ r \in \mathbbm{R}$ .  \\
    For arbitrary   {\bf hw spaces} \  $ ( X, \|..\| ) $ , \   we are able to define for \ 
	 $ \vec{x},  \vec{y} \in  X \backslash{\sf Z} $ \ with \\  
	 $ |< \vec{x} \: | \: \vec{y} >_\spadesuit| \ \leq \ \|\vec{x}\| \cdot \|\vec{y}\| $ \ \
	  an  'angle',  acording to the euclidean angle in inner product spaces. 
	  \begin{equation*}   
	   \angle_{Thy} (\vec{x}, \vec{y}) :=
                     \arccos{\frac{< \vec{x} \: | \: \vec{y} >_\spadesuit}{ \|\vec{x}\| \cdot \|\vec{y}\| } }
	  =   \arccos{ \left( \frac{1}{4} \cdot   
         \left[  \ \left\| \frac{\vec{x}}{\|\vec{x}\|}   +  \frac{\vec{y}}{\|\vec{y}\|} \right\|^{2} \  - \
         \left\| \frac{\vec{x}}{\|\vec{x}\|}   -   \frac{\vec{y}}{\|\vec{y}\|} \right\|^{2} \  \right] \right)} \ .
    \end{equation*}         
	  Then we  state that in the case of a seminormed space \ $ ( X , \|..\| ) $ , \  
	  that the triple \ \  $( X , \|..\| $,   $<.. \: | \:..>_\spadesuit  ) $  \ 
   satisfies the { \it Cauchy-Schwarz-Bunjakowsky Inequality } or { \sf CSB }inequality, \ that means  
    for all \ \ $ \vec{x}, \vec{y} \in X $ \ we have the inequality \ \
    $ |< \vec{x} \: | \: \vec{y} >_\spadesuit | \ \leq \ \|\vec{x}\| \cdot \|\vec{y}\| $ . 
    Hence in a real normed vector space \ $ ( X , \|..\| ) $ \ \ the 'Thy angle' 
     $  \angle_{Thy} (\vec{x}, \vec{y}) $ \ \  is defined for all \ $ \vec{x}, \vec{y} \neq \vec{0} $ . \\ 
   This new 'angle' has  eight nice properties, which are known from the euclidean angle in inner product spaces and      corresponds with the euclidean angle in the case that \ $ ( X , \|..\| ) $ \ already is an inner product space. \\  
     Let \  \ $ ( X , \|..\| ) $ \ be a real normed space, let dim$(X) > 1 $ , let $ \vec{x},\vec{y} \neq \vec{0}$ .  \     Then we have that
    		\begin{itemize}
     \item  \quad \ \ $\angle_{Thy}$  is a continuous surjective function 
             from \  \ $  [X \backslash \{\vec{0}\}] ^{2} $ { }  \  to \  $ [0,\pi] $ .
    \item    \quad   \    \  $\angle_{Thy}(\vec{x}, \vec{x}) = 0 $ .  
   \item     \quad  \   \ $\angle_{Thy} ( -\vec{x}, \vec{x}) = \pi$ .  
      \item   \quad   \  \  $\angle_{Thy} (\vec{x}, \vec{y}) = \angle_{Thy} (\vec{y}, \vec{x}) $ .
    \item     \quad  \ \ For all \ $ r,s > 0 $ \ ,
    we  have  \ \  $\angle_{Thy} (r \cdot \vec{x}, s \cdot \vec{y}) = \angle_{Thy} (\vec{x},\vec{y})$ .
     \item   \quad  \   \  $\angle_{Thy} ( - \vec{x}, - \vec{y}) = \angle_{Thy} (\vec{x},\vec{y})$ .              
    \item    \quad \  $ \angle_{Thy}(\vec{x}, \vec{y}) + \angle_{Thy}(-\vec{x}, \vec{y})  = \pi $ ,   
     \end{itemize} 
    which are all easy to prove. Moreover, the 'Thy angle' has  the following 
      important property, which is the main content of this paper,
    and which is not so easy to prove. \\    { \bf Theorem }    \\
     \quad  \sf  For any two  linear independent vectors \ $ \vec{x},\vec{y} $ , \  
               there is  a  decreasing    homeomorphism   
    $$   \Theta :  \mathbbm{R} \longrightarrow  (0,\pi) , \ \
                 t \mapsto \angle_{Thy}(\vec{x}, \vec{y} + t \cdot\vec{x})  \ \ . $$        
    \rm     This theorem is proved in the  section { \bf On the Existence of Polar Coordinates}. 
    We use  the help of a paper of  Charles Dimminie,  Edward Andalafte, and Raymond Freese 
     \cite{Dimminie/Andalafte/Freese1}.   \\  
   Furthermore, we work with two interesting facts from the usual two-dimensional euclidean geometry. 
   We  do not  state them here in the introduction,
    because the most difficult part is to write them down, rather than to prove them. \\
   At the end, two open questions about \ {\bf hw spaces} \ $ (X,  \|..\|) $  which have a non convex unit ball 
   are described.  One possibility to avoid these problems is 
   to replace this homogeneous weight \ $  \|..\| $ \ by another, using the convex hull of the unit ball.   
     We define a useful generalization of the introduced  \ 
    'Thy angle', maintaining all its good properties.   \\           
	 \section{General Definitions}   
    Let    \ $ X $ \ =  ( $ X , \tau$ ) \  be an  arbitrary   real topological  vector  space,  \  that  means  that 
    the real  vector space   $ X $  is  provided with a topology  $ \tau $ such that  the addition  of two vectors
    and  the  multiplication  with real numbers  are continuous.  Further let \  $\|..\|$ \ denote a 
     { \it positive functional }  on $ X$, that means that \ \  
   $\|..\|$: \ $ X \longrightarrow$  $\mathbbm{R}^{+}  \cup  \{ 0  \}$ \  is  continuous, 
   \  $\mathbbm{R}^{+} \cup \{0\}$ \ carries  the  usual  euclidean topology. 
        \\      We consider some  conditions.       \\   
       (1)  \ \         
                            For  all $  \it r \  \in \ \mathbbm{R}$ \ \rm \ and  all   \ \ 
                      \it   $\vec{x}  \in X$ \  \rm we have: \
                    \it  $\|r \cdot  \vec{x}\|$ =   $|r|$ $ \cdot  \|\vec{x}\|$ \quad  \hfill 
                                                           \rm (" absolute homogenity"),          \\
       (2) \ \    \rm  $\| \vec{x}\| = 0$ \quad  { \rm if and only if} \quad  $ \vec{x} = \vec{0}$  
                                                            \hfill  \quad \rm ("positive \ definiteness"), \\
       (3) \  \    \rm  For \ all \ \it  $\vec{x} , \vec{y} \in X$ \ \rm \    hold \  \  \ 
                   \it   $\|\vec{x}+\vec{y}\| \leq  \|\vec{x}\|  + \|\vec{y}\|$ 
                                                             \hfill \quad \rm ("triangle \ inequality"),   \\
      (4) \  \    \rm  For \ all \ \it  $\vec{x} , \vec{y} \in X$ \ \rm  \  hold \ \ \ 
                   \it   $\|\vec{x}+\vec{y}\|^{2} +    \|\vec{x}-\vec{y}\|^{2}  = 
                    2 \cdot [ \  \|\vec{x}\|^{2}  + \|\vec{y}\|^{2} \ ] $    \\  
                                \quad   $   {  } $               \hfill \quad \rm ("parallelogram identity"),               $  \begin{array}[ ]{lll}  \rm   
      \ \ \rm If \  \ \   \|..\|  \  fulfils \ (1) & \rm then  & \ \rm we \ call \ \|..\| \ \  a \  
      { \it homogeneous \ weight }  \ \rm on \ \it X,  \\  
      \ \ \rm if \ \rm \ \   \|..\|  \  fulfils \ (1) , (3) & \rm then & \  \|..\| \ \ \rm is  \ called 
       \ a \ { \it seminorm } \rm \ on \ \it X, \rm  \ \  and   \\ 
     \ \ \rm if \  \ \   \|..\|  \  fulfils \ (1), (2) \ and \ (3) & \rm then & \  \|..\| \ \  \rm is  \ called 
       \ a \ { \it norm } \rm \ on \ \it X, \rm  \ \  and   \\ 
      \ \ \rm if \  \ \   \|..\|  \  fulfils \ (1), (2) , (3) \ and \ (4) & \rm then \ the \ pair \ & \ 
      ( X,  \|..\|     )  \  \  \rm is  \ called \  an \ { \it  inner \ product \ space } . 
      
     \end{array}   $  
   
    Acording to this  cases we call                                                           
       the pair \   ( $ X , \|..\|$ ) \ a { \it homogeneously weighted vector space}  (or {\bf hw space}), 
       a  { \it seminormed vector space}, a  { \it normed vector space},  or an 
         { \it inner product space} (or {\bf IP space} ),     respectively.       \\
  	Now let \quad   $ <..|..> \ : \  X^{2} \longrightarrow   \mathbbm{R} $, \  \ let  $ <..|..> \ $ \ be continuous 
	 as  a  map  from  the product space  $ X \times X $ \ to the euclidean space  \ $ \mathbbm{R}. $
	  \qquad  We  consider some  conditions:  \\
	 $\overline{(1)}$  \ \     For  all \ \ $  \it r \  \in \ \mathbbm{R}$ \ \rm \ and  all   \ \ 
                  \it   $\vec{x}, \vec{y}   \in X $ \  \rm  hold \ \
                  \it  $ <r \cdot \vec{x} \ | \ \vec{y}>  = r \cdot <\vec{x} \ | \ \vec{y}>  
                                   $ \hfill \rm     ("homogenity"),   \\
   $\overline{(2)}$  \  \    \rm  For \ all \ \ \it  $\vec{x} , \vec{y} \in X$ \ \rm  hold \  \  \ 
                   \it  $ <\vec{x} \ | \ \vec{y}>  \  = \ <\vec{y} \ | \ \vec{x}>   $ 
                                                             \hfill \quad \rm ("symmetry"),   \\                           $\overline{(3)}$ \ \    For  all \  $  \vec{x}  \in X$ \  \rm  we have: \
                    \it  $ <\vec{x} \ | \ \vec{x}> \  \  \geq \  0  $
                                           \hfill  \quad \rm  ("positive semidefiniteness"),       \\                      $\overline{(4)}$ \ \    \rm  $  <\vec{x} \ | \ \vec{x}> \  =  0  $ \quad  { \rm if and only if} 
                                              \quad  $ \vec{x} = \vec{0}$  
                                                            \hfill  \quad \rm ("definiteness"), \\
  $\overline{(5)}$ \  \    \rm  For \ all \ \it  $\vec{x} , \vec{y} ,  \vec{z} \in X$ \ \rm   hold \quad  
           $  <\vec{x} \ | \ \vec{y}+\vec{z}> \ = \    <\vec{x} \ | \ \vec{y}>  + <\vec{x} \ | \ \vec{z}>  $   \\ 
                                \quad   $   {  } $   \hfill \quad \rm ("linearity in the second component").  
     $  \begin{array}[ ]{l}  \rm                               
       \rm If \    <.. \: | \: .. >    \  fulfils \ \overline{(1)},\overline{(2)},\overline{(3)} \ , \ \ 
    \rm then   \   \rm we \ call \  <.. \: | \: .. >   \ a \ { \it homogeneous \ product }  \ \rm on \ \it X,  \\ 
      \rm if \   <.. \: | \: .. >    \  fulfils \ \overline{(1)},\overline{(2)},\overline{(3)}, 
                                     \overline{(4)},\overline{(5)} \ , \ \ \rm then  \  \ 
         \rm  <.. \: | \: .. > \ \ is  \ called  \ an \ { \it inner \ product }  \ \rm on \ \it X.    
          \end{array}    $   \\
   Acording to this  cases we call                                                           
       the pair \   ( $ X , <.. \: | \:..>$ ) \ a  { \it homogeneous product vector space},  
             or an  { \it inner product space } (or {\bf IP space}),  respectively.       
   \begin{remark}  \rm
    We use the term  \ '{\bf IP space}' \ twice, but both definitions coincide: It is well-known that a  norm 
    is based on an inner product   if and only if the parallelogram identity holds.   
     \end{remark}                                   
   Let \   $\|..\|$ \  be denote a positive functional on $ X$. \  Then  define the two closed subsets of   $ X$:  \\
   $ {\bf S} \ := {\bf S}_{(X,\|..\|)} \ :=  \  \{ \: \vec{x} \in X \ |  \ \|\vec{x}\| = 1 \: \}$ ,
                                                                         \quad the {\it unit sphere} of $X$, \\
   $ {\bf B} \ := {\bf B}_{(X,\|..\|)} \ :=  \  \{ \: \vec{x} \in  X \ |  \ \|\vec{x}\| \leq 1 \: \} $ ,
                                                                                \quad the {\it unit ball} of $X$ .
   
    Now assume that the real vector space \ $X$ \ is provided with a positive functional \ $ \|..\|$ \ and a
   product  \  $ <.. \: | \:..> $ . \  Then the triple \  $( X , \|..\| ,  <.. \: | \:..> ) $ \ 
   satisfies the { \it Cauchy-Schwarz-Bunjakowsky Inequality } or { \sf CSB }inequality 
   $ \Longleftrightarrow $  for all \ $ \vec{x}, \vec{y} \in X $ \ we have the inequality \ \
    $ |< \vec{x} \: | \: \vec{y} >| \ \leq \ \|\vec{x}\| \cdot \|\vec{y}\| $ .  \\
    Assume  that the pair \   ( $ X , \|..\|$ ) \  is a homogeneously weighted vector space  (or {\bf hw space}). 
    Then define  for every vector \ $  \vec{v}  $ \
     with \  $ \|\vec{v}\| \neq 0 $ \  \ the vector \ \  $ { \sf sign}( \vec{v}) := 
     \frac{1}{ \| \vec{v} \| }  \cdot    \vec{v}  $, \ \       thus \ \ $ { \sf sign}( \vec{v}) $ \ \ 
   is the projection of \ \ $ \vec{v} $ \ \  into the unit sphere \ \  $ {\bf S}_{(X,\|..\|)}$ .    \\
     Let  \ $ A $ \ be an arbitrary subset  of a linear real vector space \ $ X $ . \ 
     Let $ A $ \  have the  property  that 
      for  arbitrary \ $ \vec{x} , \vec{y}  \in A $ \ and for  every  \ $ 0 \leq t \leq 1 $ \ we have \ 
	 $ t \cdot \vec{x} + (1-t) \cdot \vec{y} \in A. $  \ Such a  set $ A $ is called  { \it  convex}.
	  The unit ball \ $ {\bf B} $ \ in a seminormed space is convex because of the triangle inequality.  \\                   A convex set \ $ A $ \ in a linear topological vector space \ $ X $   =  ( $ X , \tau$ ) \ 
	 is  called { \it strictly convex } if and only if    for  all \ $ \vec{x} , \vec{y}  \in A $   \  
	 and for  every  \ $ 0 < t < 1 $ \  holds that \  \
	 $ t \cdot \vec{x} + (1-t) \cdot \vec{y} \in interior(A) $.  \\
	Let \ $A$  \ be an arbitrary subset of a real vector space $X$. Then we define the { \it convex hull } 
	of  \ $A$,  
	 $$  { \it conv(A) } := \bigcup \ \{ \sum_{i=1}^{n} t_i \cdot \vec{x_i}  \ | \ n \in \mathbbm{N} , \ t_i \in [0,1]  \ 
 \ {\rm and } \ \ \vec{x_i} \in A  \ \ { \rm for} \ i = 1, ... , n \ , \ {\rm and } \ \sum_{i=1}^{n} t_i = 1 \} \ , $$  
	 which is the smallest convex set that contains $A$.  \\
	  Let  the pair \   ( $ X , \|..\|$ ) \ be a homogeneously weighted vector space  (or {\bf hw space}),  with               the unit ball \ $ {\bf B} $ \ of $X$ . \ 
	  Let \ $ \|..\|_{| \it conv({\bf B})}$ \ be the Minkowski Functional of ${ \it conv({\bf B})}$ in $X$, \
	 that means  for all \ $\vec{x} \in X$ \ that \ \ $ \|\vec{x}\|_{| \it conv({\bf B})} \ := 
	            \ \inf \: \{ r > 0 \: | \: \frac{1}{r} \cdot \vec{x} \in { \it conv({\bf B})} \} $. \ Hence 
	            \ $ \|\vec{x}\|_{| \it conv({\bf B})} \ \leq  \  \|\vec{x}\| $.      \ \
	Note that for a 	{\bf hw space} \ $( X , \|..\| )$ , \  the pair \ 
	 $( X , \|..\|_{| \it conv({\bf B})})$ \ is a seminormed vector space. \\             
 Then we call \ $ \|..\|$ , \ or  the pair \ 
  ( $ X , \|..\|$ ), \ respectively, \   { \it normable} \ \
  if and only if  the pair \  $( X , \|..\|_{| \it conv({\bf B})})$ \ is a normed vector space.  \\
    Let \ ( $ X , \|..\|$ ) \ be a real   {\bf hw space} . \ The subset \ ${ \sf Z}$ \ of \ $X$ , \ \ 
   $ { \sf Z} :=  \{ \vec{x} \in X  \ | \ \|\vec{x}\| = 0 \} $ is  called the { \it zero-set} \ of ( $ X , \|..\|$ ) .        
   \section{On Angle Spaces} 
	In the usual euclidean plane angles are  considered for more than 2000 years. With the idea of 
	'metrics' and 'norms' others than the euclidean one  the idea came to have also 
	  orthogonality and  angles in metric and normed spaces, respectively. The first attempt to define a  concept of 
	 generalized 'angles' on metric spaces was made by Menger  \cite{Menger},p. 749 .
	 Since then a  few ideas have been developed, see the references \cite{Dimminie/Andalafte/Freese1},
	 \cite{Dimminie/Andalafte/Freese2}, \cite{J.E.Valentine/S.G.Wayment}, \cite{J.E.Valentine/C.Martin},
	 \cite{Singer}, \cite{Milicic1}, \cite{Milicic2}, \cite{Milicic3}, \cite{Gunawan}.
	 	  In this paper we focus our intention on
	 real normed spaces as a generalizitation of real inner product spaces. 
	 Let  \  $ (  X , <.. \: | \: .. > )$  be an  {\bf IP space},   
	  and let \ $ \|..\|  $ \ be 
	 the associated  norm,  $ \|\vec{x}\| := \sqrt{< \vec{x} | \vec{x} >} $, 	 then 	the  triple 
	   $( X ,  \|..\|  ,  <.. \: | \: .. >  ) $ \  fulfils the  { \sf CSB } inequality, 
	   and we have  for all \  $ \vec{x},\vec{y} \neq \vec{0} $ \ the well-known  euclidean angle   \ \   
 $\angle_{Euclid} (\vec{x}, \vec{y}) := \arccos{\frac{< \vec{x} \: | \: \vec{y} >}{ \|\vec{x}\| \cdot \|\vec{y}\|}}$ 
     \ \ with all its nice properties. 
     
  Now we want  to create a  useful definition of an 'angle space'.  Of course, if we think of angles as we
    used them in {\bf IP spaces},  we wish to get all the properties which are known from  
    these angles.     But we have to  avoid  
   extremal positions; that  means, if we demand too much of the properties of the known euclidean angle, 
    we only can expect to get 
   {\bf IP spaces} as 'angle spaces',  see \cite{Dimminie/Andalafte/Freese2}, \cite{J.E.Valentine/S.G.Wayment}.  
   On the other hand, if we request none, then we will  get a lot of  
   'angle spaces', but without any interesting characteristics. Thus we  have to find the golden mean. 
   So let's try:  
   	\begin{definition}    \rm
     Let \ ( $ X , \|..\|$ ) \ be a real   {\bf hw space} . 
      Let \  $ { \sf Z} :=  \{ \vec{x} \in X  \ | \ \|\vec{x}\| = 0 \} $ \ be the 'zero-set'. 
   We call  the triple \ $ {(  X , \|..\|,\angle_X ) } $ \  an \ { \it angle space }  if and only if 
     the following conditions ({\tt An} 1),({\tt An} 2),({\tt An} 3),({\tt An} 4),({\tt An} 5)  are satisfied.
     \newpage
   		\begin{itemize}
     \item  ({\tt An} 1) \quad $\angle_{X}$ \ \ is a continuous    
                         function from \ \ $  [X \backslash { \sf Z}] ^{2} $ \  in the interval \ $ [0,\pi] $ .
    \item    ({\tt An} 2) \quad  For all \ \ $ \vec{x} \in X \backslash { \sf Z} $ \ we have 
                                               \  $\angle_{X}(\vec{x}, \vec{x}) = 0 $ .  
   \item     ({\tt An} 3) \quad  For all \ \ $ \vec{x} \in X \backslash { \sf Z}$ \ we have 
                                 \ $\angle_{X} ( -\vec{x}, \vec{x}) = \pi$ .  
      \item  ({\tt An} 4) \quad  For all \ \ $ \vec{x},\vec{y} \in X \backslash { \sf Z}$ \ 
                           we have \  $\angle_{X} (\vec{x}, \vec{y}) = \angle_{X} (\vec{y}, \vec{x}) $ .
    \item     ({\tt An} 5) \quad  For all \ $ \vec{x},\vec{y} \in X \backslash { \sf Z}$ \ and for all \ $ r,s > 0 $ \ 
    we  have   \  $\angle_{X} (r \cdot \vec{x}, s \cdot \vec{y}) = \angle_{X} (\vec{x},\vec{y})$ .
   \end{itemize} 
    Furthermore we write down  some more properties of such conditions which seems to us 'desireable', 
   but 'not absolutely necessary'.
   	\begin{itemize}
    \item   ({\tt An} 6) \quad For all \  $ \vec{x},\vec{y}  \in  X \backslash { \sf Z}$ \  \  we  have 
            \  $\angle_{X} ( - \vec{x}, - \vec{y}) = \angle_{X} (\vec{x},\vec{y})$ .         
    \item   ({\tt An} 7) \quad   For all \  $ \vec{x},\vec{y}  \in  X \backslash { \sf Z}$ \ \   we  have 
            \  $ \angle_{X}(\vec{x}, \vec{y}) + \angle_{X}(-\vec{x}, \vec{y})  = \pi $ .
     \item   ({\tt An} 8) \quad  For all \  $ \vec{x},\vec{y} , \ \vec{x}+\vec{y} \in  X \backslash { \sf Z}$ \  \ we            have 
     \item[ ] \quad  \qquad  \  $\angle_{X}(\vec{x}, \vec{x}+\vec{y}) + 
                  \angle_{X}(\vec{x}+\vec{y},\vec{y}) =   \angle_{X}(\vec{x},\vec{y})$ .
     \item  ({\tt An} 9)  \quad                       
        For  $ \vec{x},\vec{y} , \ \vec{x}-\vec{y}  \in  X \backslash { \sf Z}$  \ \ we  have  
        \item[ ]  \quad  \qquad  \ $ \angle_{X}(\vec{x},\vec{y})  + \angle_{X}(-\vec{x}, \vec{y}-\vec{x}) +                                                   \angle_{X}(-\vec{y},\vec{x}-\vec{y})  = \pi $ .
     \item  ({\tt An} 10)   \quad      
        For all \  $ \vec{x},\vec{y} , \ \vec{x}-\vec{y}  \in  X \backslash { \sf Z}$ \  \ we  have  
       \item[ ]   \quad  \qquad  \ $\angle_{X}(\vec{y}, \vec{y}-\vec{x})
            +  \angle_{X}(\vec{x},\vec{x}-\vec{y}) =  \angle_{X}(-\vec{x},\vec{y}) $ . 
    \item ({\tt An} 11) \quad  For any two  linear independent vectors \ $ \vec{x},\vec{y}  \in X \backslash { \sf Z}$          ,     \ \   we  have  a  decreasing    
          \item[ ]       \quad  \qquad  \ \ homeomorphism \ \
    $ \Theta :  \mathbbm{R} \longrightarrow  (0,\pi) , \ \ t \mapsto \angle_{X}(\vec{x}, \vec{y} + t \cdot\vec{x})$ .          \end{itemize}      
    \end{definition}
    \begin{remark}  \rm   \label{remark3} 
   We add another demand to the above conditions. If we construct an angle \ $\angle_{Y}$ \
      for every element \
    $( Y , \|..\| ) $  \ of a class { \sf K },     and if \ \
    $ \{ ( X , \|..\| ) | ( X , \|..\| ) \ \ \rm is \ an \ { \bf IP \ space} \} \subset $ { \sf K },  \ \ 
   then  for  every {\bf IP  space}  \ $ ( Y , \|..\| ) $  
    \ \  should hold that \ \ $\angle_{Y} = \angle_{Euclid} $ .
    
      \end{remark}
    	\section{The Thy Angle}        
    Now  imagine that the real vector space \ $X$ \ is provided with a positive functional \ $ \|..\|$ \ and a
   product  \  $ <.. \: | \:..> $ .  Assume two elements $ \vec{x}, \vec{y} \in X $, 
    $ \|\vec{x}\| \cdot \|\vec{y}\| \neq 0 $ , and the property that    \  
    $ |< \vec{x} \: | \: \vec{y} >| \ \leq \ \|\vec{x}\| \cdot \|\vec{y}\| $ .  \ \ 
    Then we can define an angle between these two elements, \ \
  $ \angle (\vec{x}, \vec{y}) :=
                     \arccos{\frac{< \vec{x} \: | \: \vec{y} >}{ \|\vec{x}\| \cdot \|\vec{y}\| } }. $ \ \
    If  the triple \  $( X , \|..\| ,  <.. \: | \:..> ) $ \ 
   satisfies the { \it Cauchy-Schwarz-Bunjakowsky Inequality } or { \sf CSB }inequality, then we are able to define
   for all $ \vec{x}, \vec{y} \in X $,  $ \|\vec{x}\| \cdot \|\vec{y}\| \neq 0 $ , this angle \quad
    $ \angle_{ } (\vec{x}, \vec{y}) := 
                   \arccos{\frac{< \vec{x} \: | \: \vec{y} >}{ \|\vec{x}\| \cdot \|\vec{y}\| }} 
    \  \in [0, \pi]$ .
  
  \quad Let  the pair \ $ ( X , \|..\| \   ) $  \ \rm  be a  \ {\bf hw \ space},
      thus \ $( X , \|..\| \ ) \  \rm  fulfils \ (1)$, the absolute  homogenity. 
       \ We define a product \ \ $ < .. \: | \: ..>_\spadesuit $  \ on $ X $ .  
        \ Let for all \ $ \vec{x} , \vec{y} \in X $:  
       \begin{equation*}   
         < \vec{x} \: | \:  \vec{y} >_\spadesuit   \ \  :=  \ \             
             \begin{cases}
           0  \quad  &
       \quad \mbox{for} \quad  \|\vec{x}\| \cdot \|\vec{y}\|  =     0   \\    
        \frac{1}{4} \cdot  \|\vec{x}\| \cdot \|\vec{y}\| \cdot 
         \left[  \ \left\| \frac{\vec{x}}{\|\vec{x}\|}   +  \frac{\vec{y}}{\|\vec{y}\|} \right\|^{2} \  - \
         \left\| \frac{\vec{x}}{\|\vec{x}\|}   -   \frac{\vec{y}}{\|\vec{y}\|} \right\|^{2} \  \right]
            &  \quad \mbox{for} \quad   \|\vec{x}\| \cdot \|\vec{y}\|  \neq    0   \\            
                         \end{cases} 
     \end{equation*}
   ( Note that, in the case that \  $ ( X , \|..\| \ ) $ is already an {\bf IP space}, this definition corresponds
   with the usual definition of the inner product. ) \ \                                         
     We have  \ $ < .. \: | \: .. >_\spadesuit \ : \  X^{2} \longrightarrow    \mathbbm{R}  $ , \ \  and the properties     $\overline{(2)}$  (symmetry) \ and $\overline{(3)}$  (positive semidefiniteness) are rather trivial.
     Clearly,  $ \|\vec{x}\| =   \sqrt{ <\vec{x} \: | \: \vec{x}>_\spadesuit } \ $  for all $ \vec{x} \in X $ . 
     We  show \   $\overline{(1)}$ ,  the homogenity. 
     For a  real number  \ $ r > 0 $ \ \ it holds  that  \ \ 
               $ <r \cdot \vec{x} \ | \ \vec{y}>_\spadesuit  \ = \ r \cdot <\vec{x} \ | \ \vec{y}>_\spadesuit $ , 
               because \  $ ( X , \|..\| \   ) $  \  satisfies \ (1).    \\
     Now we prove \ \quad   
         $ <- \vec{x} \ | \ \vec{y}>_\spadesuit \ = \ - <\vec{x} \ | \ \vec{y}>_\spadesuit $ .     \ \       
    Let \  \ $ \|\vec{x}\| \cdot \|\vec{y}\|  \neq  0 $ . \  We  have 
     \begin{eqnarray*}            
       - <\vec{x} \ | \ \vec{y}>_\spadesuit     & = &   
      - \  \frac{1}{4} \cdot  \|\vec{x}\| \cdot \|\vec{y}\| \cdot 
         \left[  \  \left\| \frac{\vec{x}}{\|\vec{x}\|} + \frac{\vec{y}}{\|\vec{y}\|}  \right\|^{2} \  - \
          \left\| \frac{\vec{x}}{\|\vec{x}\|} -  \frac{\vec{y}}{\|\vec{y}\|}  \right\|^{2} \  \right]  \ \ ,  \\ 
          { \rm and }    \quad     <- \vec{x} \ | \ \vec{y}>_\spadesuit  & = &  
      \frac{1}{4} \cdot  \| -\vec{x}\| \cdot \|\vec{y}\| \cdot 
         \left[  \ \left\| \frac{ -\vec{x}}{\| -\vec{x}\|} + \frac{\vec{y}}{\|\vec{y}\|} \right\|^{2} \  - \
         \left\| \frac{ -\vec{x}}{\| -\vec{x}\|} -  \frac{\vec{y}}{\|\vec{y}\|} \right\|^{2} \  \right]    \\ 
         & = & 
        \frac{1}{4} \cdot  \|\vec{x}\| \cdot \|\vec{y}\| \cdot 
         \left[  \ \left\| \frac{\vec{y}}{\| \vec{y}\|} - \frac{\vec{x}}{\|\vec{x}\|} \right\|^{2} \  - \
         \left\| \frac{ \vec{x}}{\| \vec{x}\|} +  \frac{\vec{y}}{\|\vec{y}\|} \right\|^{2} \  \right]   \ \ ,   
   \end{eqnarray*}                
  hence \  $ <- \vec{x} \ | \ \vec{y}>_\spadesuit \ = \ - <\vec{x} \ | \ \vec{y}>_\spadesuit $ . 
  Then easily follows also for every real number $ r < 0 $ \ that \  
    $ <r \cdot \vec{x} \ | \ \vec{y}>_\spadesuit  \ = \ r \cdot <\vec{x} \ | \ \vec{y}>_\spadesuit $ , \quad  and 
     the homogenity $\overline{(1)}$ \  has been proved, hence  the pair \ 
   $ ( X , < .. | .. >_\spadesuit ) $ \ is a homogeneous product vector space.   
    
	\begin{definition}    \rm
	  For all \  {\bf hw \ spaces}  ( $ X , \|..\|$ ) \ for all  \ $ \vec{x}, \vec{y} \in X \backslash { \sf Z}  $ \  
	   ( that means $ \|\vec{x}\| \cdot \|\vec{y}\|  \neq  0 $) \ with \ 
	    $ |< \vec{x} \: | \: \vec{y} >_\spadesuit | \ \leq \ \|\vec{x}\| \cdot \|\vec{y}\| $   \ \ 
	   we define the \ { \it Thy angle } \ \ ( which is a modification of the angle  discussed in 
	   \cite{Dimminie/Andalafte/Freese1}, but there the authors  implicitly assume the  parallelogram identity.
	                                         See also \cite{Singer}. ) \\ 
	   $ \angle_{Thy} (\vec{x}, \vec{y}) :=
                     \arccos{\frac{< \vec{x} \: | \: \vec{y} >_\spadesuit}{ \|\vec{x}\| \cdot \|\vec{y}\| } }
	  =   \arccos{ \left( \frac{1}{4} \cdot   
         \left[  \ \left\| \frac{\vec{x}}{\|\vec{x}\|}   +  \frac{\vec{y}}{\|\vec{y}\|} \right\|^{2} \  - \
         \left\| \frac{\vec{x}}{\|\vec{x}\|}   -   \frac{\vec{y}}{\|\vec{y}\|} \right\|^{2} \  \right] \right)} $ .
 \end{definition}
	 \begin{proposition}  \qquad   \sf    \label{proposition1} 
	  (a)    \ \  \    If \ ( $ X , \|..\|$ ) is a real seminormed vector space, then the 	 triple \\
	      $( X , \|..\| , <.. \: | \: .. >_\spadesuit ) $ \  \  fulfils the  { \sf CSB } inequality, 
	      hence the  \ 'Thy angle' \  $ \angle_{Thy} (\vec{x}, \vec{y}) $ \\  is defined for all  \
	       $ \vec{x}, \vec{y} $ \ with \ $\|\vec{x}\| \cdot \|\vec{y}\| \neq 0 $ .                        \\ 
	  (b)  \ \ \ If \ ( $ X , \|..\|$ ) is a real seminormed vector space, then the 	 triple \ 
	  \  $(X , \|..\|,\angle_{Thy})$ \
	  fulfils  all the above demands  ({\tt An} 1), ({\tt An} 2), ({\tt An} 3), ({\tt An} 4), ({\tt An} 5) . Hence  \ \
	   ( $ X , \|..\|,\angle_{Thy} $ ) \ is  an angle space.  \\
	  (c)  \ \ \ If \ ( $ X , \|..\|$ ) is a real seminormed vector space, then the 	 triple \ 
	  \  $(X , \|..\|,\angle_{Thy})$ \	 \  fulfils \  ({\tt An} 6) \ and \ ({\tt An} 7) .   \\
	   (d)   \ \ \  \sf If \  $ (  X , <.. \: | \: .. >_{IP} )$  \ \ is an  {\bf IP space}, then 	the  
	    triple \  $( X , \|..\| ,  <.. \: | \: .. >_{IP}  ) $ \ \  fulfils the  { \sf CSB } inequality 
	   and we have   \   for all \  $ \vec{x},\vec{y} \neq \vec{0} $ \  that \ \   
	  $ \angle_{Thy} (\vec{x}, \vec{y})  \  =  \ \angle_{Euclid} (\vec{x}, \vec{y})  $ .   \\
	  (e)\ \ \ If \ ( $ X , \|..\|$ ) is a real normed vector space, then the 	 triple \ 
	  \  $(X , \|..\|,\angle_{Thy})$ \ generally does not  fulfil \ ({\tt An} 8), ({\tt An} 9), ({\tt An} 10) .  \\
   (f)\ \ \ If \ ( $ X , \|..\|$ ) is a real seminormed vector space, then the 	 triple \ \
   $ (X,\|..\|,\angle_{Thy})$ \ generally does not  fulfil \ 
   ({\tt An} 8), ({\tt An} 9), ({\tt An} 10), ({\tt An} 11) .            \\ 
	  \end{proposition}
	   \begin{proof}
	 	  (a) \quad If \ ( $ X , \|..\|$ ) is a real seminormed vector space, then because of the triangle inequality 
	  and \ \ $ \left\| \frac{\vec{x}}{\|\vec{x}\|} \right\| = 1 $  \ \  
	   we get that \ \  $ \left| < \vec{x} \: | \:  \vec{y} >_\spadesuit \right| =  
	  \left|  \frac{1}{4} \cdot  \|\vec{x}\| \cdot \|\vec{y}\| \cdot 
         \left[  \ \left\| \frac{\vec{x}}{\|\vec{x}\|}   +  \frac{\vec{y}}{\|\vec{y}\|} \right\|^{2} \  - \
         \left\| \frac{\vec{x}}{\|\vec{x}\|}   -   \frac{\vec{y}}{\|\vec{y}\|} \right\|^{2} \  \right] \right| \\
	 \leq  \  \frac{1}{4} \cdot  \|\vec{x}\| \cdot \|\vec{y}\| \cdot 
	 \max \left\{  \left\| \frac{\vec{x}}{\|\vec{x}\|} + \frac{\vec{y}}{\|\vec{y}\|} \right\|^{2},
	               \left\| \frac{\vec{x}}{\|\vec{x}\|} - \frac{\vec{y}}{\|\vec{y}\|} \right\|^{2}  \right\} 
	 \leq     \frac{1}{4} \cdot  \|\vec{x}\| \cdot \|\vec{y}\| \cdot 2^{2}  
	   \  =  \   \|\vec{x}\| \cdot \|\vec{y}\| $  .  \\
	 (b) \quad  Rather trivial if you use (a) and the fact that  \ $\|..\|$ \ is homogeneous.     \\
	 (c) \quad  ({\tt An} 6) is trivial because \  $< \: | \: >_\spadesuit$ \ is homogeneous,   and ({\tt An} 7) is easy      if you   know that \ $ \arccos(r) + \arccos(-r) = \pi $ .  \\
	 (d) \quad If \ ( $ (  X , <.. \: | \: .. >_{IP} )$  is an  {\bf IP space} with the associated 
	 norm  $ \|\vec{x}\| := \sqrt{ < \vec{x} \: | \vec{x} \: >_{IP}} $ , 
	 then,  because of \ \ $  <..\:|\:.. >_\spadesuit  \ = \ <.. \: | \: .. >_{IP} $ , \  we have  \  \
	   $ \angle_{Thy} (\vec{x}, \vec{y})  \  =  \ \angle_{Euclid} (\vec{x}, \vec{y})  $ .  \\
	 (e) \quad We need counterexamples.  Recall the pairs  ( $ \mathbbm{R}^{2} , \|..\|_{p}$) ,  with 
	 the  { \it Hölder weights } $ \|..\|_{p}$, \ $ p > 0$ , \ we define that \ 
	  $ \|(x_1,x_2)\|_{p} := \  \sqrt[p]{|x_1|^{p}+|x_2|^{p}}$  . \ The pairs  ( $ \mathbbm{R}^{2} , \|..\|_{p}$) 
	   are normed spaces if and only if $ p \geq 1$ . 
	 For \ $ p=2$ \ we get the usual euclidean norm. So let us take, for instance, \ $ p = 1 $ , because it is easy 
	 to calculate with.  \\   
	 Let \ $ \vec{x} := (1,0), \ \vec{y} := (0,1) $ , both vectors have  the \ $ \|..\|_{1}$-norm 1 . 
	                                                     \quad  Then  we  have
	 \begin{eqnarray*}
	   \angle_{Thy} (\vec{x}, \vec{y})  
	    &  = &  \arccos{ \left( \frac{1}{4} \cdot   
         \left[  \ \left\| \frac{\vec{x}}{\|\vec{x}\|_1}   +  \frac{\vec{y}}{\|\vec{y}\|_1} \right\|_1^{2} \  - \
         \left\| \frac{\vec{x}}{\|\vec{x}\|_1}   -   \frac{\vec{y}}{\|\vec{y}\|_1} \right\|_1^{2} \  \right] \right)} \\
     & =  &  \arccos{ \left( \frac{1}{4} \cdot   
         \left[  \ \left\| (1,0)  + (0,1) \right\|_1^{2} \  - \
         \left\|  (1,0)  - (0,1) \right\|_1^{2} \  \right] \right)}    \\   
     &  =  &   \arccos{ \left( \frac{1}{4} \cdot   
         \left[  \ 4 \  - \ 4 \  \right] \right)} \ = \ \arccos (0) \ = \ \pi / 2 \ = \ 90 \deg   .
    \end{eqnarray*}   
    \begin{eqnarray*}
	  { \rm  And }  \quad \quad   \angle_{Thy} (\vec{x},\vec{x} + \vec{y})  
	     & =  &  \arccos{ \left( \frac{1}{4} \cdot   
         \left[  \ \left\| (1,0)  +  \frac{1}{2} \cdot  (1,1) \right\|_1^{2} \  - \
         \left\|  (1,0)  - \frac{1}{2} \cdot  (1,1) \right\|_1^{2} \  \right] \right)}    \\   
      & =  &  \arccos{ \left( \frac{1}{4} \cdot  \left[  \ \left( 2 \right)^{2} \  - \ 
        \left( 1 \right)^{2}  \  \right] \right)}  =  \arccos \left(\frac{3}{4}\right) \ \approx  \ 41.41 \deg   .
    \end{eqnarray*} 
    With similar calculations, we get \ \ 
    $  \angle_{Thy} (\vec{x} + \vec{y},\vec{y})  \ = \  \arccos \left(\frac{3}{4}\right) $ ,  \quad
    hence \\  $ \ \angle_{Thy} (\vec{x},\vec{x} + \vec{y}) \ + \ \angle_{Thy} (\vec{x} + \vec{y},\vec{y})  
     \ \neq \ \angle_{Thy} (\vec{x},\vec{y})  $ , \  and    that contradicts \ ({\tt An} 8).   \\
    The condition \ ({\tt An} 9) \ means that the sum of the inner angles of a triangle is $\pi$.    \\
    We can use the same example of the normed space \  ($ \mathbbm{R}^{2} , \|..\|_{1}$) \ with  unit vectors
     $ \vec{x} := (1,0)$ , and \ $\vec{y} := (0,1) $. Again we get \\ 
     $\angle_{Thy} (\vec{x},\vec{y}) = \pi/2 ,  \ \ 
     \angle_{Thy} (-\vec{x},\vec{y}-\vec{x}) =  \angle_{Thy} (-\vec{y},\vec{x}-\vec{y}) 
                                                                =  \arccos(\frac{3}{4}) $ , \ \ hence \\
     $ \angle_{Thy} (\vec{x},\vec{y}) +  
     \angle_{Thy} (-\vec{x},\vec{y}-\vec{x}) +  \angle_{Thy} (-\vec{y},\vec{x}-\vec{y}) < \pi $ ,
     hence \ ({\tt An} 9) \ is not fulfiled.  \\
  For the condition ({\tt An} 10)  we use the same space and the same vectors 
  $ \vec{x} := (1,0)$ , and \ $\vec{y} := (0,1)$.
  \  We get \   $\angle_{Thy} (-\vec{x},\vec{y}) = \pi/2 ,  \ \ 
     \angle_{Thy} (\vec{y},\vec{y}-\vec{x}) =  \angle_{Thy} (\vec{x},\vec{x}-\vec{y}) 
                                                                =  \arccos(\frac{3}{4}) $, 
      \ \ hence \ ({\tt An} 10) \ is not fulfiled. \\              
	(f) \quad We use the same example as in (e) to prove that \ ({\tt An} 8), ({\tt An} 9), ({\tt An} 10) \  generally is     not fulfiled.  Or we can change it, so that it is no more  a normed space. Take   the pair 
	  ( $ \mathbbm{R}^{3} , \|..\|_{\widehat{1}}$) , with   \ 
	$ \|..\|_{\widehat{1}} (x,y,z) := |x| + |y| $ . \ Obviously, it is a seminormed,  but not a normed space,
	   with  does not fulfil   \ ({\tt An} 8), ({\tt An} 9), ({\tt An} 10).  \\
	  Here is a further example.  
	 \ \ Let \ $  (  X , \|..\| ) := ( \mathbbm{R}^{2},\|..\| ) $ \ be the seminormed space with the seminorm \
	   \ $ \|(x_1,x_2)\| :=  |x_1|$. \ Hence \  $ { \sf Z} = \{ (0,x_2) \: | \: x_2 \in \mathbbm{R} \}$ \
	    is the zero-set.  \  We get only two angles,
	   for all  \ $ \vec{x},\vec{y} \in \mathbbm{R}^{2}\backslash { \sf Z}$ \ hold that \  
	   $ \angle_{Thy} (\vec{x}, \vec{y})  \in \{ 0, \pi \} $.  \ 
	   Then ({\tt An} 11) is not satisfied: Take \ \ $ \vec{x} := (1,0),  \ \vec{y} := (1,1) $ , \ then for all 
	   $ t \in \mathbbm{R}\backslash \{-1\} $ \ we have that for \ $ t > -1 $ , \ 
	   $ \angle_{Thy} (\vec{x}, \vec{y}+ t \cdot \vec{x} ) = 0 $, \ and  for \ $ t < -1 $ , \  \ 
	     $ \angle_{Thy} (\vec{x}, \vec{y}+ t \cdot \vec{x} ) = \pi $ .  The calculations for this  
	     are easy.
	  	\end{proof} 
	 Another interesting non-trivial example is the following.  \\
	    Let \ \ $  (  X , \|..\| ) := ( \mathbbm{R}^{2},\|..\| ) $ \ \ be a {\bf hw space} with the unit sphere \ \ 
	  $ {\bf S}  :=  \{ \: \vec{x} \in  \mathbbm{R}^{2} \ |  \ \|\vec{x}\| = 1 \: \} 
	  \ := \ \{ (x_1,x_2) \in 	\mathbbm{R}^{2} \ | \ |x_2| \cdot |x_1| = 1 \ \} $. \
	 Hence \  $ { \sf Z} = \{ (0,x_2) \: | \: x_2 \in \mathbbm{R} \} \cup \{ (x_1,0) \: | \: x_1 \in \mathbbm{R} \}$ \
	   \  is the zero-set.  \ \ This space   fulfils the  { \sf CSB } inequality, and \ \ 
	   	  \  $( \mathbbm{R}^{2}, \|..\|,\angle_{Thy})$ \  satisfies   ({\tt An} 1),  
	   	  {({\tt An} 2)}, ({\tt An} 3), ({\tt An} 4), ({\tt An} 5), ({\tt An} 6), ({\tt An} 7) .
	   	  Hence  \ \  ( $ \mathbbm{R}^{2} , \|..\|,\angle_{Thy} $ ) \ is 
	   	      an angle space, which is not a seminormed space ( the unit ball is not convex).  \  We have  \ 
	   $ \angle_{Thy} (\vec{x}, \vec{y})  \in \{ 0, \pi/2, \pi \} $ \ 	      
	    for all  \ $ \vec{x},\vec{y} \in \mathbbm{R}^{2}\backslash { \sf Z}$ .  
	 \newpage  
	\section{On the Existence of Polar Coordinates}   
	\begin{theorem}  \sf   \label{main Theorem} 
	   Assume that \  ( $ X , \|..\|$ ) is a real normed vector space. Then the Thy  angle   $\angle_{Thy} $                \ \  satisfies \ \  ({\tt An} 11) . In other words, for a linear independent subset 
	    \ $ \{\vec{x},\vec{y} \} \subset X $ \ we get a   decreasing  homeomorphism \ 
    $ \Theta :  \mathbbm{R} \longrightarrow  (0,\pi) , \ \ t \mapsto \angle_{Thy}(\vec{x}, \vec{y} + t                      \cdot\vec{x})$ . 
		\end{theorem} 
	The above theorem is 	the main result of this paper.    	
	Before we start the lengthy proof ( finished on page 17 ) we formulate two comments.
	\begin{remark}   \rm
	The theorem remains  true when we formulate it more general, but less beautiful:  \\
	 '  Assume that \  ( $ X , \|..\|$ ) is a real seminormed vector space.
	  Assume that a linear independent subset \ $ \{\vec{x},\vec{y} \} \subset X $ \  generates the 
	  two-dimensional subspace $ U  \subset X $, assume \   $ { \sf Z} \cap U = \{\vec{0}\} $ . 
          Then we  have  a  decreasing  homeomorphism \ 
    $ \Theta :  \mathbbm{R} \longrightarrow  (0,\pi) , \ \ t \mapsto \angle_{Thy}(\vec{x}, \vec{y} + t                      \cdot\vec{x})$ ' .
   \end{remark} 
	\begin{corollary}    \sf
	 With the above theorem we can describe elements of a two-dimensional real normed vector space \
	   $(  X , \|..\| )$ \  by polar coordinates. 
	 If we fix a basis \ $ \{ \vec{b_1},\vec{b_2} \} $ , \ 
	  then  every \ $ \vec{x} = r_1 \cdot \vec{b_1} + r_2 \cdot \vec{b_2} \: \in X $ \ 
	  is uniquely defined  by its norm  \ $ \|\vec{x}\| $ \ and its angle  \quad  
	 $  \angle_X (\vec{x}) := \angle_{Thy} (\vec{x}, \vec{b_1}) $ \ \ if and only if \ \ $ r_2 > 0 $  \ \  and  \ \
	 $  \angle_X (\vec{x}) := -\angle_{Thy} (\vec{x}, \vec{b_1}) $ \ \ if and only if \ \ $ r_2 < 0 $ .
	 This concept easily can be extended to finite dimensional real normed vector spaces. 
	\end{corollary}	
	Now we start the proof of the above theorem. \
	The reader should have a copy of the remarkable paper
	of   Charles Dimminie,  Edward Andalafte, and Raymond Freese  \cite{Dimminie/Andalafte/Freese1}, 
	because we need some propositions from that paper, which we write down without the proofs.  \\ 
	The central idea of the proof is, assuming that the map $ \Theta $ is not injective contradicts the 
	convexity of the unit ball \ $ {\bf B} $ \ of \ $ ( X , \|..\|) $ . 
	\begin{proof} \  ( of { \bf Theorem \ref{main Theorem}} )  \\ 
	 Assume that \  ( $ X , \|..\|$ ) is a real normed vector space. ( Hence $ { \sf Z} = \{\vec{0}\} )$ .
	  Assume that a linear independent subset \ $ \{\vec{x},\vec{y} \} \subset X $ \  generates the 
	  two-dimensional subspace $ U $. We consider   the map \ $\Theta :  \mathbbm{R} \longrightarrow  [0,\pi] $ .
	  For  convenience, we define  some  abbreviations. \  Let $ h_+ , h_- :  \mathbbm{R} \longrightarrow  [0,2] ,
	 \quad    h_+(t) := \left\| \frac{\vec{x}}{\|\vec{x}\|} +  
	                              \frac{\vec{y} + t \cdot\vec{x}}{\|\vec{y} + t \cdot\vec{x}\|} \right\| , \ \     
	          h_-(t) := \left\| \frac{\vec{x}}{\|\vec{x}\|} -  
	                              \frac{\vec{y} + t \cdot\vec{x}}{\|\vec{y} + t \cdot\vec{x}\|} \right\|$ .
	                              \qquad We have
	  \begin{eqnarray*}
	   \Theta(t)  :=  \angle_{Thy}(\vec{x}, \vec{y} + t\cdot\vec{x}) & = &                                                      
	       \arccos{ \left( \frac{1}{4} \cdot       \left[  \ \left\| \frac{\vec{x}}{\|\vec{x}\|}   +                              \frac{\vec{y} + t \cdot\vec{x}}{\|\vec{y} + t \cdot\vec{x}\|} \right\|^{2} \  - \
       \left\| \frac{\vec{x}}{\|\vec{x}\|}   -   
       \frac{\vec{y} + t \cdot\vec{x}}{\|\vec{y} + t \cdot\vec{x}\|} \right\|^{2} \  \right] \right)}  \\ 
      & = &  \arccos{ \left( \frac{1}{4} \cdot      
	        \left[  \left[    h_+(t) \right]^{2} \  - \ \left[ h_-(t) \right]^{2} \ \right]  \right)}  \ \  ,       
       \end{eqnarray*}
    and let  \quad $  \Theta(-\infty) := \pi , \quad  \Theta(+\infty) := 0 $ .    
	  \begin{lemma}  \sf
	  We have that  \\  
	    $ \lim_{t \rightarrow -\infty}   h_+(t) =  \lim_{t \rightarrow +\infty}   h_-(t) = 0 $ , \quad \sf and \quad
	   $ \lim_{t \rightarrow +\infty}   h_+(t) =  \lim_{t \rightarrow -\infty}   h_-(t) = 2 $ .
	  \end{lemma} 
	  \begin{proof}
	  See  \cite{Singer},p.38 ,  \ or \  \cite{Dimminie/Andalafte/Freese1},p.199 .  
	  \end{proof}
	   \begin{lemma}  \sf
	  We have that  the map \ 
	     $\Theta :  \mathbbm{R} \cup \{-\infty,+\infty\} \longrightarrow  [0,\pi] $ \ 
	     is     continuous and surjective.
	    \end{lemma} 
	  \begin{proof} \ \
	   By the previous lemma  \quad $ \lim_{t \rightarrow -\infty}  \Theta(t) = \pi $ \quad and 
	                          \quad      $ \lim_{t \rightarrow +\infty}  \Theta(t) = 0 $ , \quad
	  and the norm $ \|..\| $ \ is continuous, hence \ $ \Theta  $ \ is continuous and
	  the image of \ $ \Theta $ \ is \ $ [0,\pi] $ .     
	      \end{proof}
	 We  still have to prove the injectivity of $\Theta $ , the difficult part of the proof.  In the  following 
	 we need the  notion of a strictly convex set, that is a subset \ $ A $ \  of a linear topological vector space
	  such that \ $ t \cdot \vec{x} + (1-t) \cdot \vec{y} \in interior(A) $ \
	  for  arbitrary \ $ \vec{x} , \vec{y}  \in A $ \ and for  every  \ $ 0 < t < 1 $ .	  
	 \begin{lemma}  \sf   \ \ 
	 The maps \ \ $ h_+ , h_- :  \mathbbm{R}  \rightarrow [0,2] $ are monotone increasing, respectively 
	                                              decreasing.  \\
	 In the case that the unit ball \ $ {\bf B} $ of $ (X, \|..\|)$ \  is strictly convex, then 
	     $ h_+ , h_- $ \ are strictly monotone  increasing, respectively decreasing.    
    \end{lemma} 
      \begin{proof}  \ 
	  See  the tricky proof in  \cite{Dimminie/Andalafte/Freese1},p.201, theorem 2.4 and p.202, theorem 2.5 . 
	  The authors  only dealt with the map \ $ h_- $ . Note that they only prove that  $ h_- $ is monotone. The case 
	  of a strictly convex unit ball  \ $ {\bf B} $  is not explicitly written down.  You have to read 
	   both proofs of the theorems 2.4 and 2.5  attentively.  
	  \end{proof}
  \begin{remark}  \rm
   If we have  a unit ball \ $ {\bf B} $ of $ (X, \|..\|)$ \ that  is not strictly  convex, then 
	     $ h_+ , h_- $ \ are generally not strictly monotone.  This is shown by the example of  the normed space \
	 $ (\mathbbm{R}^{2},\|..\|_{\infty})$ \ with the norm  $ \|(x_1,x_2)\|_{\infty} := \max \{ |x_1|,|x_2| \} $ .  
	 Choose \ $ \vec{x} := (1,0)$ , and \ $ \vec{y} := (0,1)$ ,  then
	 the interval in which \ $h_+$ \ is constant \ $( =1 ) $ \ is \ $ [-1,0] $ , \ while   
	  \ $h_-$ \ is constant \ $( =1 ) $ \ in \ $ [0,1] $. This example shows that the intervals where 
	    \ $h_+$ \ and \ $h_-$, \ respectively, are constant may intersect in one point. We are just proving the fact 
	    that both intervals do not intersect in an  interval with nonempty interior.  
    \end{remark}  
  \begin{corollary}    \sf
   In the case that the unit ball \ $ {\bf B} $ of a real normed space $ (X, \|..\|)$ \  is strictly convex,
   then   \    ({\tt An} 11) is  satisfied.
  \end{corollary} 
    \begin{proof}
    We had defined  \quad 
     $   \Theta(t)  =  \arccos{ \left( \frac{1}{4} \cdot      
	   \left[  \left[    h_+(t) \right]^{2} \  - \ \left[ h_-(t) \right]^{2} \ \right]  \right)}   $ ,  
	   \ \  and because  $ {\bf B} $   is strictly convex it follows by the last lemma that 
	   \ $ h_+, h_- $ \ are strictly monotone increasing, respectively decreasing.   
	 \end{proof}
	  Now we  want to prove that \ $  \Theta  $ \  remains to be strictly monotone decreasing,
	  even  if the unit ball is 	 not strictly convex.  Because \ $ h_-, h_+ $ \ are monotone, 
	  \ $  \Theta  $ \ is always monotone decreasing. We have to prove that the monotony is strict.  \\
	  We prefer a direct proof. \ \ 
	 Assume that \ $ -\infty < t_1 \leq t_2 < +\infty $ \ , and  \ $  \Theta(t_1) =  \Theta(t_2) $  . \\
	 The case \ $ t_1 = t_2 $ \ is possible.  We will show that this is the only possible case. 
	 So  let us assume that   \ $ -\infty < t_1 < t_2 < +\infty $ . \  Now we  hunt for contradictions. 
	 Because of  \  $  \Theta(t_1) =  \Theta(t_2) $ \ and because  \ $ h_+$ \ is monotone increasing and \
	  \ $ h_-$ \ is monotone decreasing,  we have that \
	 \ $ h_+(t_1) =  h_+(t_2) $ \  and \  \ $ h_-(t_1) = h_-(t_2) $ ,  \ and for all \ $ t $ \ in the interval 
	 $ [ t_1, t_2 ] $ \ both $ \Theta $ and $ h_+ $ and $ h_- $ remain constant.  
	 
	 For the rest of the proof we calculate in coordinates of the basis   \
	  $ \{ \frac{\vec{x}}{\|\vec{x}\|}, \frac{\vec{y}}{\|\vec{y}\|} \} \ = \
	  \{ { \sf sign}( \vec{x}), { \sf sign}( \vec{y}) \} $.   \\
	  Hence  \quad $ { \sf sign}( \vec{x})  $ = 
	    $  \left(  \begin{array}{c}  1 \\ 0   \end{array}  \right)  $ ,  \ \ 
	     $ { \sf sign}( \vec{y})   =   \left(  \begin{array}{c} 0 \\ 1    \end{array}  \right)  $  \ \
	      $ \vec{v} :=    \left(  \begin{array}{c}       v_{1} \\ v_{2}     \end{array}  \right)  $ 
	      := $  { \sf sign}( \vec{y} + t_1 \cdot\vec{x})  $ ,   \ \  
	    $   \vec{w} :=      \left(  \begin{array}{c}   w_{1} \\ w_{2}  \end{array}  \right)  $ 
	    :=   $  { \sf sign}( \vec{y} + t_2 \cdot\vec{x})  $ .   \ \   
	We have  \ $  \vec{v} =      \left(  \begin{array}{c}   v_{1} \\ v_{2}   \end{array}  \right)  $ 
	= $   \frac{\vec{y} + t_1 \cdot\vec{x}}{\|\vec{y} + t_1 \cdot\vec{x}\|} $                                      
	= $   \frac{t_1 \cdot \|\vec{x}\|}{\|\vec{y} + t_1 \cdot\vec{x}\|} \cdot  \frac{\vec{x}}{\|\vec{x}\|} +             
	      \frac{ \|\vec{y}\|}{\|\vec{y} + t_1 \cdot\vec{x}\|} \cdot  \frac{\vec{y}}{\|\vec{y}\|} $ ,
	  hence \ $ v_2 =   \frac{ \|\vec{y}\|}{\|\vec{y} + t_1 \cdot\vec{x}\|} > 0 $ .                                           For the same reason we get \ \ $  w_2 = \frac{ \|\vec{y}\|}{\|\vec{y} + t_2 \cdot\vec{x}\|} > 0 $ . 
	  \begin{lemma}    \qquad  \sf  We have \ \   $ \vec{v} \neq \vec{w} $ .
	  \end{lemma}
	   \begin{proof}   \quad
	  We assume the opposite \ $ \vec{v} = \vec{w} $ . \ Because of  \
	   $ v_2 = w_2$ , \ we have \\ $ v_2 =  \frac{ \|\vec{y}\|}{\|\vec{y} + t_1 \cdot\vec{x}\|} = 
	                              w_2 =  \frac{ \|\vec{y}\|}{\|\vec{y} + t_2 \cdot\vec{x}\|} $ , \ hence \ 
	      $ \|\vec{y} + t_1 \cdot\vec{x}\| = \|\vec{y} + t_2 \cdot\vec{x}\| $ .  \ Then follows with \  $ v_1 = w_1$ 
	 \ that \ \ $ t_1 \cdot \|\vec{x}\| = t_2 \cdot \|\vec{x}\| $ , \ hence \ $ t_1 = t_2 $ .                           
	   \end{proof}
	  Because of \ \ $ h_+(t_1) =  h_+(t_2) $ \ \ and \  \ $ h_-(t_1) = h_-(t_2) $ , \  we have that \ \ 
	   $   h_+(t_1) =  \left\| \left(  \begin{array}{c} 
	                v_1 + 1 \\ v_2 
	               \end{array}  \right)  \right\|  $ =  $ h_+(t_2)  = 
	                  \left\| \left(  \begin{array}{c} 
	                w_1 + 1 \\ w_2 
	               \end{array}  \right) \right\|  $ 
	   \ \ and \ \    
	   $  h_-(t_1) =  \left\| \left(  \begin{array}{c} 
	                v_1 - 1 \\ v_2 
	               \end{array}  \right)  \right\|  $ =  $  h_-(t_2)  = 
	                  \left\| \left(  \begin{array}{c} 
	                w_1 - 1 \\ w_2 
	               \end{array}  \right) \right\|  $  .   \\     \\
	   For further investigations we  must distinguish a few cases. 
	   \begin{itemize}                              
	    \item  \ Case ({ \bf A}):   $  v_1 < w_1  $ , \quad with the subcases 
	      \begin{itemize}  
	             \item  \ Case ({ \bf A1}):  \quad  $ v_1 < w_1 < -1 $ \quad or \quad $ +1 < v_1 < w_1   $ .
	             \item  \ Case ({ \bf A2}):  \quad $   v_1 < w_1 $ \quad and \quad  $ \{v_1, w_1\} \cap [-1,+1] \neq                                                                                       \emptyset $
	     \end{itemize}
       \item \ Case ({ \bf B}):   $  v_1 > w_1  $ , \quad with the corresponding subcases \ ({ \bf B1}),({ \bf B2}) .
       \item \ Case ({ \bf C}):   $ v_1 = w_1 $ , \quad with the subcases 
               \begin{itemize}  
	             \item  \ Case ({ \bf C1}):  \quad $   v_1 = w_1  \in \ \{-1,+1\} $ ,
	             \item  \ Case ({ \bf C2}):  \quad $   v_1 = w_1 \in (-1,+1) $ ,
	             \item  \ Case ({ \bf C3}):  \quad  $ v_1 = w_1 < -1 $ \quad or \quad $ v_1 = w_1 > +1 $ .
	             \end{itemize}
    \end{itemize}
  Let's start with  the  easy \quad Case ({ \bf C1}): \  \ Assume, for instance,  \ \ $  v_1 = w_1  = 1 $ . \ \           Because of   \ $ h_-(t_1) =  h_-(t_2) $ \  it follows  \                               %
     $ \left\| \left(  \begin{array}{c} 
	                 0 \\ v_2  
	               \end{array}  \right) \right\|   =     \left\| \left(  \begin{array}{c} 
	                0 \\ w_2 
	               \end{array}  \right) \right\|  $ ,  \ hence \ $ v_2 = w_2 $ . \ 
	               This contradicts \  $ \vec{v} \neq \vec{w} $ . \ \  Or more detailed, \ \   $ v_2 = w_2 $ 
	               \ means \ $ \| \vec{y} + t_1 \cdot \vec{x} \| = \| \vec{y} + t_2 \cdot \vec{x} \| $ , \ hence \ 
	                $ \frac{t_1 \cdot \|\vec{x}\|}{\|\vec{y} + t_1 \cdot\vec{x}\|} = v_1 = 1 = 
	                w_1 = \frac{t_2 \cdot \|\vec{x}\|}{\|\vec{y} + t_2 \cdot\vec{x}\|}   $ , \   which is 
	                  only possible for  \    $  t_1 = t_2  $ ,  and there is a contradiction.  \\  \\
	  For all further cases, note that we can replace \ 
	                          $ t_1, t_2 $ \  by \ $ \widetilde{t_1}, \widetilde{t_2}, $ \ with \ 
	                             $ t_1 <  \widetilde{t_1} < \widetilde{t_2} < t_2 $ \  and \ 
	$ \vec{v} :=  \left(  \begin{array}{c}   v_1 \\ v_2         \end{array}  \right)  := 
	   \frac{\vec{y} + \widetilde{t_1} \cdot\vec{x}}{\|\vec{y} + \widetilde{t_1} \cdot \vec{x}\|} $,   
	$  \vec{w} :=  \left(  \begin{array}{c}   w_1 \\ w_2         \end{array}  \right)  := 
	   \frac{\vec{y} + \widetilde{t_2} \cdot\vec{x}}{\|\vec{y} + \widetilde{t_2} \cdot \vec{x}\|}  $  \ \  
	         to make sure that  all seven unit vectors  \ \ 
	  \begin{displaymath}    \left(  \begin{array}{c}    0 \\ 1    \end{array}  \right) ,
	   \ \  \left(  \begin{array}{c}   v_1 \\ v_2    \end{array}  \right), \  
	    \left(  \begin{array}{c}   w_1 \\ w_2         \end{array}  \right),  \   
	            { \sf sign}\left(  \begin{array}{c} v_1 - 1 \\ v_2     \end{array}  \right), \  
	            { \sf sign}\left(  \begin{array}{c} v_1 + 1 \\ v_2     \end{array}  \right),  \
	            { \sf sign}\left(  \begin{array}{c} w_1 - 1 \\ w_2     \end{array}  \right), \ 
	             { \sf sign}\left(  \begin{array}{c} w_1 + 1 \\ w_2     \end{array}  \right)                                 \end{displaymath} \quad  are  distinct.  \\    \\
	  Now we deal with the even more easier  \quad Case ({ \bf C2}): \qquad Thus \ \ $ -1 <  v_1 = w_1  < 1 $ .                          Because of \ \ $ \vec{v} \neq \vec{w} $ \ \  we have \ \ $ v_2 \neq w_2 $ , \ \  for instance \ \ 
	  $ 0 < v_2 < w_2 $ . \    Assume \ \ $ 0 \leq  v_1 = w_1  < 1 $ .  Hence, by convexity of 
	  the unit ball \ $ {\bf B} $ \ of \ $ ( X , \|..\|) $ ,  \ \ and $ \vec{w} , 
	   \left( \begin{array}{c} 1 \\ 0 \end{array} \right)  \in  {\bf S} $ \ , the straight line between both points 
	   is in   \ $ {\bf B} $ . \  But  then $ \vec{v} $  would be in the interior of  $ {\bf B} $ .
	   This  is impossible,  because \ $ \|\vec{v}\| = 1 $ , (see \ Picture 1).  \\  
	   \setlength{\parindent}{0mm}
    \begin{figure}[h]
    \centering
     \setlength{\unitlength}{1cm}
    \begin{picture}(6,2)                 \put(3.3,1.5){  $ \vec{w}$ }   \put(3.3,0.8){  $ \vec{v}$ }  
       \put(0,0){\vector(1,0){6.0}}     \put(3,-0.3){\vector(0,1){3.0}}  
       \put(0.75,-0.4){$-1$}     \put(1.2,-0.1){\line(0,1){0.2}}     \put(2.55,1.65){ $  1 $}  
       \put(4.5,-0,35){ $  1 $}                                      \put(2.9,1.8){\line(1,0){0.2}} 
       \put(3.5,2.5){\line(1,-2){1.5}}    \put(3.65,1.53){  $ \times $ }    \put(3.65,0.8){  $ \times $ }  
         \put(6.15,-0,3){$   { \sf sign}( \vec{x})  $}     \put(1.75,2.3){$   { \sf sign}( \vec{y})  $}                          \put(2.5,-1.0){Picture 1}      
      \end{picture}
    \end{figure}	
              \\    \\    \\   \\    \\    \\  
	   Now we need a break.  Let us turn  our attention to 'Grecian Geometry', with this phrase we mean
	    elementary  geometry on the two-dimensional  euclidean plane. 
	    We  describe and prove two propositions about  'projections', which will be needed for the proof
	    of the theorem.      \\     
	  For a better understanding of the following propositions it is reasonable to look  
	  at the enclosed drawings. 
	  \begin{proposition}   \sf      \label{proposition A}
    Let us take \ $ \mathbbm{R}^{2} = \{ ({\tt x}|{\tt y}) | \ {\tt x},{\tt y} \in \mathbbm{R} \} $ ,
     \ the two-dimensional  euclidean plane, with the  horizontal {\tt x}-axis and  the vertical  {\tt y}-axis.
     Consider the two parallel lines \
      $  G_S: {\tt y} = {\tt x} - 1 \ $ { \rm and} \ $ G_T: {\tt y} = {\tt x} + 1 \ $ . 
     Assume a third  straight   line \ $ L $ ,  \ not parallel to \ $  G_S, G_T $, respectively, 
      \ 
        with the property     that  \ $ L $ \ does not meet the origin $(0|0)$.  
     The  intersection of \ $ L $ \ with \ $ G_S $ \ is called \ $ S = ({\tt x}_S|{\tt x}_S-1) $, \ \ 
     and the intersection of \ $ L $ \ with \ $ G_T $ \ is called \  $ T = ({\tt x}_T|{\tt x}_T+1) $ . \\
     Then  there is an unique point \ 
      $  P_{hor} = ({\tt x}_{hor}|{\tt y}_{hor})$ \ on \ $L$ , \  
        such that the three points  \\   $ (0|0), ({\tt x}_{hor}+1 | {\tt y}_{hor})$, \ and \ $ S $ \ 
        are collinear, \quad and such that the three points   \ \ 
           $ (0|0), ({\tt x}_{hor}-1 | {\tt y}_{hor})$, \ and  $ T $ \ 
        are collinear. ( See the enclosed Picture 2 ).
      \end{proposition}
      \begin{proof} In the case of \ \ $ {\tt x}_S \neq {\tt x}_T $, we have an equation \ 
      $ L: {\tt y} = m_L \cdot {\tt x} + b_L  , \ m_L, b_L  \in \mathbbm{R}, \ b_L \neq 0 , \ m_L \neq 1 $, \ 
        and then take  
       \begin{eqnarray*}   
        P_{hor} \ = \ ( \ {\tt x}_{hor}|{\tt y}_{hor}) & = & 
         \left( \frac{m_L - {b_L}^{2}}{b_L \cdot ( m_L -1 ) }  
                        \ | \  \frac{ {m_L}^{2} - {b_L}^{2}}{b_L \cdot ( m_L -1 ) } \right)    \\
     & = &  \left( \ \frac{{\tt x}_T \cdot ( b_L + 1 ) + 1}{b_L}  
     \ | \  m_L \cdot \frac{{\tt x}_T \cdot ( b_L + 1 ) + 1}{b_L} + b_L \ \right)  \\ 
      & = &  \left( \ \frac{{\tt x}_S \cdot ( b_L - 1 ) + 1}{b_L}
     \ | \   m_L \cdot \frac{{\tt x}_S \cdot ( b_L - 1 ) + 1}{b_L} + b_L \ \right) \ \ ,  
    \end{eqnarray*}   
   $  { \rm with } \quad  m_L \ = \ \frac{ 2 + {\tt x}_T - {\tt x}_S } {{\tt x}_T - {\tt x}_S} , \quad 
   { \rm and } \quad  b_L \ = \ {\tt x}_S \cdot ( 1 - m_L ) - 1  \ = \ {\tt x}_T \cdot ( 1 - m_L ) + 1 $ , \ 
   \ and some elementary calculations confirm that indeed all three formulas of \ $ P_{hor} $ \  
   yield the same values for \ $ {\tt x}_{hor} $ \ and \ ${\tt y}_{hor}$ , \ and that  \ $ P_{hor} $ \ fulfils 
   the demanded properties.   \ \ ( See Picture 2 ) .  \\
   In the case of \ \ $ {\tt x}_S = {\tt x}_T $ ,   \ \ 
    we have an equation \ \ $ L: \ {\tt x} = a_L := {\tt x}_S = {\tt x}_T $ , \ and  we get 
    \begin{eqnarray*}   
        P_{hor} \ = \ ( \ {\tt x}_{hor}|{\tt y}_{hor}) \ = \ \left( a_L \ | \ a_L - \frac{1}{a_L} \right)  \ ,
        \ \ \ { \rm  ( see \ Picture \ 7 ) . }
    \end{eqnarray*}  
       \end{proof}  
     \setlength{\parindent}{0mm}
    \begin{figure}[ht]
    \centering
     \setlength{\unitlength}{1cm} 
    \begin{picture}(6,6)       \put(2.0,-1.5){Picture 2}   
      \thinlines  
      \put(6.15,-0,1){$   {\tt x}   $}     \put(-0.3,5.8){$  {\tt y} $}     \put(-0.3,4.8){$5$}  
       \put(-6.0,0){\vector(1,0){12.0}}     \put(0,-3){\vector(0,1){9}}  
       \put(1.0,-0.3){$1$}     \put(3.1,3.9){$G_T$}  
                  \put(5.1,3.9){$G_S$}   \put(6.1,-1.5){$L: {\tt y}=-{\tt x}+5$} 
       \put(-1.55,0,15){ $ -1 $}    \put(-0.2,1.0){$1$}    \put(0.1,-1.2){$-1$}  
         \put(-4.0,-3.0){\line(1,1){7.0}}  \put(-2.0,-3.0){\line(1,1){7.0}}   \put(-1.0,6.0){\line(1,-1){7.0}} 
            \put(2.85,1.55){$S$}     \put(1.84,3.2){$T$}    
               \put(-1.0,-1.5){\line(2,3){4.5}}  \put(-4.0,-2.66666666){\line(3,2){8.5}} 
          \put(4.5,2.9){ ${\tt y}= \frac{2}{3}{\tt x}$}    
         \put(3.5,4.9){${\tt y}= \frac{3}{2}{\tt x} $}     
                    \put(2.2,2.65){$P_{hor}$ }    \put(2.45,2.30){$\times$ }  
                            
      \put(0.8,6.6){ Here: \  \ $ m_L = -1, \ b_L = 5$ . } 
      \put(0.8,5.9){ We get $S=(3|2) $ \  and \ $T=(2|3) $ . \ Hence \ $P_{hor}=(\frac{13}{5}|\frac{12}{5}) $ . }  
       \linethickness{0.5mm} \put(1.6,2.4){\line(1,0){2.0}}           
      \end{picture}
    \end{figure}	            
     \newpage 
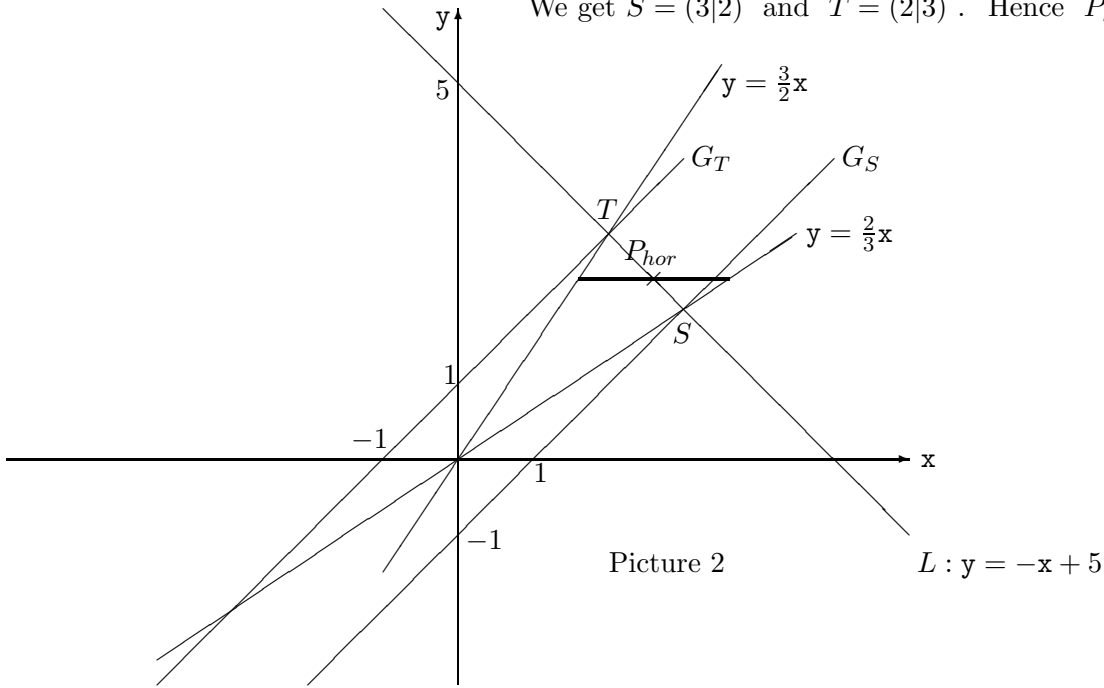
   \begin{proposition}   \sf    \label{proposition B}
        Let us again take \ $ \mathbbm{R}^{2} $ \ 
     with the  horizontal {\tt x}-axis and  the vertical  {\tt y}-axis.
  Consider the straight line \ $ G  $ \ with the equation \   \ 
  ${\tt y} = m \cdot {\tt x} + 1 , \ m \in \mathbbm{R} \backslash \{-1,+1\} $ .  \ 
    Let us choose  an arbitrary point \ 
      $ ( \widehat{{\tt x}} | \widehat{{\tt y}} )$ \ on $ G  , \ \widehat{{\tt y}} \neq 0 $ . \ 
    Take two points \  $ S := ( \widehat{{\tt x}}-1 | \widehat{{\tt y}} )$ \ and 
                    \  $ T := ( \widehat{{\tt x}}+1 | \widehat{{\tt y}} )$. \
    We call \ $ \overline{S}$ \ the projection of \ $ S $ \ on the line $ G $ , \  and  $ \overline{T}$ \ 
    the projection   of \ $ T $ \ on the line $ G $. \ (That means that the three points \
      $ (0|0), S , \overline{S}$, \ and the three points \ $ (0|0), T , \overline{T} $, respectively,
      \ are collinear, \ $ \overline{S}, \overline{T}  \in G $.)  The four points \ 
      $  \overline{S}, \overline{T},  -\overline{S}, -\overline{T} $ \ are the corners of a parallelogram.
      We claim that the intersection \ $ \nu $ \ of the line that connects \ $ \overline{T} $ \ and \
       $  -\overline{S}$ \   with the    horizontal  {\tt x}-axis has the value \ $1$, \   and  the intersection of the
    line that connects \ $ -\overline{T} $ \  and \ $  \overline{S}$ \   with the  {\tt x}-axis is \ $ -1$ .  
     Moreover, the two points \ \ $ \overline{S}, \overline{T} $ \ are on the same  side of the horizontal
     {\tt x}-axis \ if and only if \ $ -1 < m < 1$ .             
     \end{proposition} \quad   (See Picture 3). 
           \\ 
	   \setlength{\parindent}{0mm}
    \begin{figure}[ht]
    \centering
     \setlength{\unitlength}{0.9cm}
    \begin{picture}(8,6)     \thinlines  
      \put(9.15,-0,1){$   {\tt x}   $}     \put(-0.3,5.8){$  {\tt y} $}    
       \put(-6.0,0){\vector(1,0){15.0}}     \put(0,-3){\vector(0,1){9}}  
          \put(-2.1,1.5){$ ( \widehat{{\tt x}} | \widehat{{\tt y}} )$ }    \put(0.15,3.9){$ 1 $ }  
                                \put(-4.4,-0.4){$ -\nu=-1 $ }              \put(3.8,-0.4){$ \nu=1 $ }  
           \put(1.2,6.1){$G: {\tt y} = 2\cdot{\tt x}+1$}   
            \put(-3.5,-3.0){\line(1,2){5.0}}   \put(-5.333333333333,1.3333333333333){\line(1,0){8.0}} 
           \put(-2.6,-1.7){$ \overline{T} $}   \put(-1.5,0.5){$ \overline{S} $}   \put(1.35,-0.88){$ -\overline{S} $}    
         \put(2.8,1.0){$ T = -\overline{T}$}   \put(-5.6,1){$ S $}     
        \put(-5.64,1.41){\line(4,-1){10.0}}   \put(-3.666666,-1.5333333){\line(5,1){10.0}} 
              \put(-5.0,-0.2){\line(5,1){10.0}}   \put(-4,-2){\line(2,1){10.0}}  
      \put(4.8,6.3){Here \ $ G: {\tt y}=2\cdot{\tt x}+1  $ , hence $ m = 2$ , }    
      \put(4.8,5.6){and \ $ (\widehat{{\tt x}}  | \widehat{{\tt y}}) := (-\frac{1}{3}|\frac{1}{3}) $ . } 
     \put(4.8,4.9){We have \ \  $ \overline{S }=( -\frac{4}{9}|\frac{1}{9} ) ,    \quad 
                           \overline{T}=(-\frac{2}{3}|-\frac{1}{3})$ .}  
      \put(4.8,4.2){Hence we  get \ \ $ \nu = 1 $ . }   
    %
    %
    %
            
      \end{picture}
    \end{figure}	       
    \quad  \\  \\ \\  \begin{center}  Picture 3    \end{center} 
       \quad  \\   
    \begin{proof} With elementary calculations, we have that  \\ 
      $$  \overline{S} \ =  \ \frac{1}{1 + m } \cdot 
    ( \widehat{{\tt x}} - 1 \: | \: \widehat{{\tt y}} ) \ \ \quad  { \rm and }  \  \quad 
      \overline{T} \ =  \ \frac{1}{1 - m } \cdot 
    ( \widehat{{\tt x}} + 1 \ | \ \widehat{{\tt y}} ) \ . $$ 
    Some more calculations yield the formula  \quad 
    $ {\tt y} = \frac{ m \cdot \widehat{{\tt x}} + 1 } {\widehat{{\tt x}} + m }
        \cdot [ \; {\tt x} - 1 \; ]  $  \quad   for the straight line that intersects \ 
       $  \overline{T}  \   \rm  and \  -\overline{S} \ , $  and finally we get  \quad  $ \nu = 1 $ . \\ 
       And \  $ \overline{S}, \overline{T} $ \ are on the same  side of the horizontal
     {\tt x}-axis \ if and only if \  their second components have the same signs, hence if and only if \ \ 
     $ (1-m) \cdot (1+m) > 0$ . 
        \end{proof} 
      Before  we  can return to our main purpose, ( that is to prove  { \bf Theorem  \ref{main Theorem}}) , we have        to  mention one  easy fact about unit balls in the $ \mathbbm{R}^{2}$ .  
  \begin{lemma}   \sf    \label{lemma A}
  Assume that \ $\mathbbm{R}^{2}$ \ is provided with a seminorm \ $ \left\|..\right\|$ \ 
   with the unit ball \  $ {\bf B} $ .  \quad  Define two closed sets \  $ Set1 $ \ and \  $ Set2$ ,  \\
  $ Set1 := \{  ({\tt x}|{\tt y}) \in \mathbbm{R}^{2} \quad | \ \ \ \ {\tt x}+1 \geq {\tt y} \geq  {\tt x}-1 \} $ , \ \
  $ Set2 := \{  ({\tt x}|{\tt y}) \in \mathbbm{R}^{2} \quad | \  -{\tt x}+1 \geq {\tt y} \geq  -{\tt x}-1 \} $  . 
   \newpage
  Assume now that we have the four unit vectors \ \ $ (1|0), (-1|0), (0|1), (0|-1) $ , \ that means \
  \ $ \left\| (1|0) \right\| =  \left\| (-1|0) \right\| =  \left\| (0|1) \right\| =  \left\| (0|-1) \right\| 
               \  = \ 1 $  . \\ Then we have that \qquad  ${\bf B} \ \subset \ Set1 \cup Set2 $ . 
  \end{lemma}                    
  \begin{proof}
  Instead of a proof we prefer to show a picture and we remark that \ $ {\bf B} $ \ has to be convex.  
   \qquad \qquad    (  \rm  See  Picture 4 ).       \end{proof} 
      	   \setlength{\parindent}{0mm}
    \begin{figure}[ht]
    \centering
     \setlength{\unitlength}{1cm}
    \begin{picture}(6,6)     \thinlines  
      \put(6.15,-0,1){$   {\tt x}   $}     \put(-0.3,5.8){$  {\tt y} $}    
       \put(-6.0,0){\vector(1,0){12.0}}     \put(0,-3){\vector(0,1){9}}

         \put(-4.0,-3.0){\line(1,1){8.0}}  \put(-2.0,-3.0){\line(1,1){7.0}}  
          \put(-5.0,4.0){\line(1,-1){7.0}}   \put(-4.0,5.0){\line(1,-1){8.0}}  \put(-2.2,1.7){ $ Set2$ } 
         \put(1.7,-2.2){ $ Set2$ }  \put(-2.3,-2.0){ $ Set1$ }     \put(1.6,2.0){ $ Set1$ } 
           \put(5.1,4){${\tt y}={\tt x}-1$}    \put(4.1,5){${\tt y}={\tt x}+1$}                     
           \put(-6.1,4.1){${\tt y}=-{\tt x}-1$}    \put(-5.1,5.1){${\tt y}=-{\tt x}+1$} 
        \put(-1.45,-0.45){ $-1$ }   \put(0.8,-0.4){ $ 1 $ }    \put(0.04,-1.1){ $-1$ }   \put(0.1,0.9){ $1$ } 
             
      \end{picture}
    \end{figure}	       
   \  \\ \\ \\ \\ \\  \\ \begin{center}  Picture 4    \end{center} 
   \quad  \\ 
  The time has come to return to the proof of \ { \bf Theorem \ref{main Theorem}} ,
   but we still need some general preparations. Recall that we had two unit vectors \\
   $  \vec{v} =      \left(  \begin{array}{c}   v_{1} \\ v_{2}   \end{array}  \right)  $ 
	  = $   \frac{\vec{y} + t_1 \cdot\vec{x}}{\|\vec{y} + t_1 \cdot\vec{x}\|}   \ \  \neq \ \  
	   \vec{w} =      \left(  \begin{array}{c}   w_{1} \\ w_{2}   \end{array}  \right)  $ 
	  = $   \frac{\vec{y} + t_2 \cdot\vec{x}}{\|\vec{y} + t_2 \cdot\vec{x}\|} $ , \ \ with $ t_1 < t_2 $ ,
	  \ \ and \ \  $ v_2, w_2 > 0 $ . \\   
   Because they are different, they uniquely determine a straight line $L$ that connects them.  If we define 
   two lines \ $L_-$ \ and \  $L_+$ \ such that \  $L_-$ \ connects \ 
   $ \vec{v} - (1|0) = ( v_1 - 1 | v_2 )$ \  and  \ $ \vec{w} - (1|0) = ( w_1 -1 | w_2 )$ , \ \ and  such that \ \
   \ $L_+$ \ connects \ $ ( v_1 + 1 | v_2 )$ \  and  \ $ ( w_1 + 1 | w_2 )$ , it is trivial 
   that all the three  lines \ $ L$ , $ L_-$ ,  and  $ L_+$ \ are parallel, \ ( see  Picture 6 ) .  
   \begin{lemma}   \sf
   Assume an arbitrary \ {\bf hw space} \ $ (Y,\left\|..\right\|) $ \ with a homogeneous weight $ \left\|..\right\|$. 
   \ Let \ $ \vec{a}, \vec{b} \in Y$ \ be linear independent, \ and let \ $ \left\|\vec{a}\right\| 
   = \|\vec{b}\|  >  0$ . \ Consider the  two-dimensional subspace of \ $Y$ , generated by the vectors \ 
   $ \vec{a}, \vec{b}$ , which is isomorphic  to the vector space \ $\mathbbm{R}^{2}$.    \\
   Then the line that connects \ $ \vec{a}$ \ and \ $ \vec{b} $ \ is parallel to the line that
    connects \ $  { \sf sign}(\vec{a})$ \ and \ $  { \sf sign}(\vec{b}) $. 
   \end{lemma}
   \begin{proof} \quad   Trivial by the intercept theorem and the fact that \ $ \left\|..\right\|$ \ is 
      homogeneous. \\  (See Picture 5) . 
   \end{proof}   
     \setlength{\parindent}{0mm}
    \begin{figure}[ht]
    \centering
     \setlength{\unitlength}{1cm}
    \begin{picture}(4,4)     \thinlines  
      \put(6.15,-0,1){$   {\tt x}   $}     \put(-0.3,4.8){$  {\tt y} $}    
       \put(-6.0,0){\vector(1,0){12.0}}     \put(0,-2){\vector(0,1){7}}  
      
       \put(-3.0,-2.0){\line(1,2){4.0}}  \put(-2.0,-2.0){\line(1,2){3.0}}  
       
            \put(-1.7,1.2){$\vec{a}$}     \put(-1.1,2.2){$\vec{b}$}   
             \put(-1.13333333,0.15){${\sf sign}(\vec{a})$}      \put(-0.4,1.0){${\sf sign}(\vec{b})$} 
      \put(0,0){\vector(-1,2){0.5}}   \put(0,0){\vector(-1,2){1.0}}  
             \put(0,0){\vector(-3,2){1.5}}       \put(0,0){\vector(-3,2){0.75}}
         \end{picture}
    \end{figure}
       { $ $ } \\  \\ \\ \\  { $ $ } \\  
         \begin{center}  Picture 5    \end{center}   { $ $ } \\	\\ 
     Further we define two straight lines \  $L_{-,{ \sf sign}}$ \ and \  $L_{+,{ \sf sign}}$ , \ 
    such that \  $L_{-,{ \sf sign}}$ \  connects \ 
    $ { \sf sign}( v_1 - 1 | v_2 )$ \  and  \ $ { \sf sign}( w_1 -1 | w_2 )$, \ \ and such that \ \ $L_{+,{ \sf sign}}$ 
    \  connects \ $ { \sf sign}( v_1 + 1 | v_2 )$ \  and  \ $ { \sf sign}( w_1 + 1 | w_2 )$ .  
    \begin{lemma}     \sf
    We claim that  all the five lines \ \ $ L, \ L_-, \ L_+, \ L_{-,{ \sf sign}},$ 
     \ and \ $ L_{+,{ \sf sign}} $ \ are parallel.   
   \end{lemma} 
   \begin{proof} \quad   \rm 
   On page 9                          
   we described the  conditions that  we have  \ \ 
	   $   h_+(t_1) =  \left\| \left(  \begin{array}{c} 
	                v_1 + 1 \\ v_2 
	               \end{array}  \right)  \right\|  $ =  $ h_+(t_2)  = 
	                  \left\| \left(  \begin{array}{c} 
	                w_1 + 1 \\ w_2 
	               \end{array}  \right) \right\|  $ 
	   \ \ and \ \    
	   $  h_-(t_1) =  \left\| \left(  \begin{array}{c} 
	                v_1 - 1 \\ v_2 
	               \end{array}  \right)  \right\|  $ =  $  h_-(t_2)  = 
	                  \left\| \left(  \begin{array}{c} 
	                w_1 - 1 \\ w_2 
	               \end{array}  \right) \right\|  $ , \\  and  all norms are greater than \ 0 , 
	    \   hence together with the previous lemma the claim is true.  
    \end{proof} 
    \begin{lemma}      \sf     \label{lemma F}
    Recall that on page 9 \                   
     we had defined six  unit vectors
     \begin{displaymath}  
	   \ \  \left(  \begin{array}{c}   v_1 \\ v_2    \end{array}  \right), \  
	    \left(  \begin{array}{c}   w_1 \\ w_2         \end{array}  \right),  \   
	            { \sf sign}\left(  \begin{array}{c} v_1 - 1 \\ v_2     \end{array}  \right), \  
	            { \sf sign}\left(  \begin{array}{c} v_1 + 1 \\ v_2     \end{array}  \right),  \
	            { \sf sign}\left(  \begin{array}{c} w_1 - 1 \\ w_2     \end{array}  \right), \ 
	             { \sf sign}\left(  \begin{array}{c} w_1 + 1 \\ w_2     \end{array}  \right) \ ,                               \end{displaymath} 
	  and we mentioned on page 10 \ that, without loss of generality, all these vectors  are distinct, and they are 
	   different from \ $(0|1)$ .    \quad     Then all the six  points are collinear.  \\ 
    In the case that these six vectors are on different  sides of the vertical  {\tt y}-axis, even all seven
    vectors
     \begin{displaymath}    \left(  \begin{array}{c}    0 \\ 1    \end{array}  \right) ,
	   \ \  \left(  \begin{array}{c}   v_1 \\ v_2    \end{array}  \right), \  
	    \left(  \begin{array}{c}   w_1 \\ w_2         \end{array}  \right),  \   
	            { \sf sign}\left(  \begin{array}{c} v_1 - 1 \\ v_2     \end{array}  \right), \  
	            { \sf sign}\left(  \begin{array}{c} v_1 + 1 \\ v_2     \end{array}  \right),  \
	            { \sf sign}\left(  \begin{array}{c} w_1 - 1 \\ w_2     \end{array}  \right), \ 
	             { \sf sign}\left(  \begin{array}{c} w_1 + 1 \\ w_2     \end{array}  \right)                                 \end{displaymath}  \\   are collinear.  
    \end{lemma}
  \begin{proof}
   All these seven unit vectors  have a positive second component, thus they are above the horizontal {\tt x}-axis.
   By the previous lemma, the three lines \ \ $ L, \  L_{-,{ \sf sign}},$ \ and \ $ L_{+,{ \sf sign}} $ \ 
   are parallel. \ \ Let us consider, for instance, \ $ L $ \ and \ $ L_{-,{ \sf sign}}$ . \
   $L$ \ meets \ $ \vec{v}$ \ and \ $ \vec{v}$ , and \
    \  $L_{-,{ \sf sign}}$ \ meets \ $ { \sf sign}( v_1 - 1 | v_2 )$ \  and  \ $ { \sf sign}( w_1 -1 | w_2 )$ . \
   Because of the convexity of the unit ball \ $ {\bf B} $ , the case
   \ $ L \neq   L_{-,{ \sf sign}}$ \ is not possible.   Hence all  six points that generate these
   three lines have to be on the   same straight line  \ $L \ = \ L_{-,{ \sf sign}} \ = \ L_{+,{ \sf sign}} $ . \  
   If in addition these six vectors are on different  sides of the  vertical  {\tt y}-axis, the unit vector
    \ $ (0|1) $ \ has to be on the same  line \ $L$, too .  \\ (See Picture 6). 
  \end{proof} 
    \setlength{\parindent}{0mm}
    \begin{figure}[ht]
    \centering
     \setlength{\unitlength}{1cm}
    \begin{picture}(4,4)     \thinlines  
      \put(6.15,-0,1){$   {\tt x}   $}     \put(-0.3,4.3){$  {\tt y} $}    
       \put(-6.0,0){\vector(1,0){12.0}}     \put(0,-1.5){\vector(0,1){6.1}}  
       
     \put(-6.0,1.0){\line(3,1){10.0}}    \put(-7.0,1.66666666){\line(3,1){6.0}}  \put(-2.0,1.33333){\line(3,1){6.0}} 
        \put(-3.6,1.74){$\times$ }      \put(-2,2.305){$\times$ }
        \put(-3.5,1.5){$\vec{v}$ }      \put(-2,2.0){$\vec{w}$ }   
         \put(-0.5,1.5){$\vec{v}+(1|0)$ }      \put(+1,2.0){$\vec{w}+(1|0)$ } 
           \put(-6.5,1.5){$\vec{v}-(1|0)$ }      \put(-6,2.55){$\vec{w}-(1|0)$ } 
                          \put(-6.43,1.85){\line(1,0){6.0}}       \put(-4.85,2.4){\line(1,0){6.0}} 
          \put(4.1,4.3){$L$ }   \put(-1.3,3.76666){$ L_{-}$ }    \put(3.8,3.4666666){$ L_{+}$ } 
             \put(-0.8,3.1){$(0|1)$}    \put(2.9,-0.45){1}   
             \put(-3.4,-0.45){$-1$}  
      \put(3,-0.1){\line(0,1){0.2}}      \put(-3,-0.1){\line(0,1){0.2}}   \put(0,0){\vector(0,1){3.0}}     
          \end{picture}
    \end{figure}	         \          \\    \\   \\  
    \begin{center}  Picture 6    \end{center}   
         { $ $ }   \\  { $ $ } 
      Now we have collected all the  facts we will need in the following. Recall that we wanted to prove 
   { \bf Theorem  \ref{main Theorem}  },  and  that we already have proved  Case ({ \bf C1}) and Case ({ \bf C2}). 
   Still there are  missing the cases  { \bf C3}, { \bf A1}, { \bf A2}, { \bf B1}, and { \bf B2} . \\       	    
	 Recall that we calculate  with the basis \ 
	  $ \{ \frac{\vec{x}}{\|\vec{x}\|}, \frac{\vec{y}}{\|\vec{y}\|} \} \ = \
	  \{ { \sf sign}( \vec{x}), { \sf sign}( \vec{y}) \} $ , \ \ and that we  have \\  
	  $  \vec{v} =      \left(  \begin{array}{c}   v_{1} \\ v_{2}   \end{array}  \right)  $ 
	  = $   \frac{\vec{y} + t_1 \cdot\vec{x}}{\|\vec{y} + t_1 \cdot\vec{x}\|}   \ \  \neq \ \  
	   \vec{w} =      \left(  \begin{array}{c}   w_{1} \\ w_{2}   \end{array}  \right)  $ 
	  = $   \frac{\vec{y} + t_2 \cdot\vec{x}}{\|\vec{y} + t_2 \cdot\vec{x}\|} $ , \ \ with $ t_1 < t_2 $ ,
	  \ \ and \ \  $ v_2, w_2 > 0 $ . \\    Because of \ $ \vec{v} \neq \vec{w}$ \ there is an unique 
	  straight line $L$ that connects both points,  with an equation  \ 
 $L: {\tt y}=m_L \cdot {\tt x}+ b_L \ , m_L, b_L \in \mathbbm{R}$ \ \ or \ \ $ L: {\tt x} = a_L $ , \
     ( if \  $ v_{1} =  w_{1} =: a_L $)  .  \\   \\
	 \ Case ({ \bf C3}):  \quad  $ v_1 = w_1 < -1 $ \ \ or \ \ $ +1 < v_1 = w_1  $ , for instance, assume \
	 $ 1 < v_1 = w_1 $ . Thus  \ $L$ \ is vertical, \ \ $ L: {\tt x} = a_L :=  v_1 = w_1 $ .  \\ 
	 Thus we deduce with proposition  \ref{proposition A} that there is an unique point \ $P_{hor}$ , 
	  \ ( see Picture 7 ) ,   
	   \begin{eqnarray*}   
        P_{hor} \ = \ ( \ {\tt x}_{hor}|{\tt y}_{hor}) \ = \ \left( a_L \ | \ a_L - \frac{1}{a_L} \right)  \ ,
   \end{eqnarray*}   
	 such that the three points \ $ (0|0), (a_L - 1 | a_L - \frac{1}{a_L} )$ , and \ $ T := ( a_L | a_L + 1 )$ 
	 \ are collinear, 
	 and such that the three points \ $ (0|0), (a_L + 1 | a_L - \frac{1}{a_L} )$ , and \ $ S := ( a_L | a_L - 1 )$ \
	 are collinear.  Because of \ $a_L > 1 $, \ $P_{hor}$ \ is between $S$ and $T$.  \
	 By the last lemma  \ref{lemma F}, all six unit vectors 
	      $$ \vec{v}, \ \vec{w}, \ { \sf sign}( \vec{v}+(1|0)), 
	           \ { \sf sign}( \vec{w}+(1|0)), \ { \sf sign}( \vec{v}-(1|0)), \ { \sf sign}( \vec{w}-(1|0)) $$
	 are located on the same line \ $L$, and they have to be in the set \ $ Set1 $
	  ( by lemma  \ref{lemma A}),    
	 and above the \  {\tt x}-axis.  Now we have to use proposition \ref{proposition A} . \  
	  Assume that \ $ \vec{v} := P_{hor}$ . \\
	  If we have a \ $\vec{w}$ \  from the point \  $  P_{hor} $ \ on \ $L$ \ upward,
	  the point \ ${ \sf sign}( \vec{w}-(1|0)) $ \ will be above \ $T$ , \ hence not in the set \ $ Set1 $ .
	     If we have a \ $\vec{w}$ \  from \  $  P_{hor}$ \ on \ $L$ \
	     downward, \ ${ \sf sign}( \vec{w}+(1|0)) $ \ will be below \ $S$, hence not in  $ Set1 $ .    
	  Thus the only possibility is \  $  P_{hor} \ = \ \vec{v} \ = \ \vec{w}$, \ 
	  ( and \  $ \ { \sf sign}( \vec{v}+(1|0))= { \sf sign}( \vec{w}+(1|0))= S$ \ and \  
	 $\ { \sf sign}( \vec{v}-(1|0))= { \sf sign}( \vec{w}-(1|0))= T$ ) , and we find a contradiction, 
	    and \ Case ({ \bf C3}) \ is discussed.     \newpage     
	   \setlength{\parindent}{0mm}
    \begin{figure}[ht]
    \centering
     \setlength{\unitlength}{1cm}
    \begin{picture}(6,6)     \thinlines  
      \put(6.15,-0,1){$   {\tt x}   $}     \put(-0.3,5.8){$  {\tt y} $}    
       \put(-6.0,0){\vector(1,0){12.0}}     \put(0,-3){\vector(0,1){9}}  
       \put(1.0,-0.3){$1$}     \put(3.1,3.9){${\tt y}= {\tt x}+1$}  
                  \put(5.1,3.9){${\tt y}= {\tt x}-1$}   \put(2.1,-1.5){$L: {\tt x}=a_L$} 
       \put(-1.45,0,15){ $-1$}    \put(-0.36,1.0){$1$}    \put(0.1,-1.2){$-1$}  
         \put(-4.0,-3.0){\line(1,1){7.0}}  \put(-2.0,-3.0){\line(1,1){7.0}}   \put(2.0,-2.0){\line(0,1){7.0}} 
            \put(2.05,0.75){$S$}     \put(1.8,3.0){$T$}    
               \put(-1.0,-1.5){\line(2,3){4.5}}  \put(-3.0,-1.5){\line(2,1){8.5}} 
          \put(-2.4,-2){$ Set1 $ }     
                    \put(1.85,1.42){$\times$ }     \put(1.97,1.65){$P_{hor}$ } 
                            
      
      \linethickness{0.5mm} \put(1.0,1.5){\line(1,0){2.0}}           
      \end{picture}
    \end{figure}	         { $ $ } \\  \\  \\   
	    \quad  \\   \begin{center}  Picture 7    \end{center}    
	   \quad { $ $ }  \\ 
	 The next  Case ({ \bf A1}) will be proved in a similar way.  \\
	 We  have \quad  $ v_1 < w_1 < -1 $ \quad or \quad $ +1 < v_1 < w_1 $ , assume \ $ 1 < v_1 < w_1$ .
	 The two different points \ $\vec{v}, \vec{w}$ \  are in the set \ $ Set1 $, and they  define an 
	 unique straight line \
	  $ L: {\tt y} = m_L \cdot {\tt x} + b_L \ , \ m_L,b_L \in  \mathbbm{R}$. If $ \ m_L = 0$ \ we have \
	   \ $ b_L \geq 1$ \ ( otherwise, if \ $ 0 < b_L < 1$ , it contradicts the convexity of the  unit ball  $ {\bf B}$, 
	   note the unit vectors \ $ (0|1), \vec{v}, \vec{w}$ ). \
	  If $ \ m_L \neq 0$ , \ we have a zero \ $a_L := -b_L/ m_L$ \ of \ $L$.  It must be both \ 
	  $ 1 \leq |a_L| $ \ and \ $ 1 \leq |b_L| $ . \ In all other cases,
	   namely \ $ -1 < a_L < 1 $ \ or \ $ -1 < b_L <  1 $ , we have a contradiction to the convexity of  $ {\bf B}$. 
	 \  ( Note the six unit vectors \ $ (1|0), (-1|0), (0|1), (0|-1), \vec{v}, \vec{w}$ ).    
	  Hence, \ $ 1 \leq |b_L| $ \ and, if \ $ m_L \neq 0$ , \ $ 1 \leq |a_L| $ . \ 
	   \begin{lemma}  \sf   Furthermore, we have \ \ $ m_L \neq 1$ . 
	   \end{lemma}  
	  \begin{proof} \ \ Assume \  $ m_L = 1$ . \
	  \ If \ $ \vec{v}, \vec{w} \in  interior(Set1) $, \
	 it contradicts the convexity of  \ $ {\bf B}$, (note \ $ (0|1), (1|0), \vec{v}, \vec{w} $). \ Or, if \
	 $\vec{v}, \vec{w} \in   L: {\tt y} = {\tt x} + 1 $ \ or \ $\vec{v}, \vec{w} \in   L: {\tt y} = {\tt x} - 1 $ , \
	 with \ $  \left\|\vec{v}+(1|0)\right\| =  \left\|\vec{w}+(1|0)\right\| $ \ or \ 
	 $ \left\|\vec{v}-(1|0)\right\| =  \left\|\vec{w}-(1|0)\right\| $ \ always follows \
	  $ \vec{v} =  \vec{w} $ .   
	 \end{proof}
	 We call \ $ S := ({\tt x}_S|{\tt x}_S-1)$ \ the intersection of \ $L$ \ and \ ${\tt y} = {\tt x} - 1$ , \ and \ 
	$T := ({\tt x}_T|{\tt x}_T+1) $ \ the intersection  of \ $L$ \ and \ ${\tt y} = {\tt x} + 1$, $ S, T \in Set1 $ .   
	      Then, by proposition \ref{proposition A},          
	        we have an unique point \
	     $  P_{hor} = ({\tt x}_{hor}|{\tt y}_{hor})$ \ on \ $L$ , \  
        such that the three points  \quad   $ (0|0), ({\tt x}_{hor}+1 | {\tt y}_{hor})$, and \ $ S $ \ 
        are collinear, \quad and such that the three points \ $ (0|0), ({\tt x}_{hor}-1 | {\tt y}_{hor})$, and  $ T $ \ 
        are collinear, ( see again Picture 2 for an example) , \ and \  $  P_{hor}$ \ has the representation 
	   
     $$   P_{hor} \ = \ ( \ {\tt x}_{hor}|{\tt y}_{hor})  =  \left( \frac{m_L - {b_L}^{2}}{b_L \cdot ( m_L -1 ) }  
                            \ | \  \frac{ {m_L}^{2} - {b_L}^{2}}{b_L \cdot ( m_L -1 ) } \right) .   $$    
	  By lemma  \ref{lemma A}   
	  and lemma  \ref{lemma F} ,   
	  all six unit vectors $$ \vec{v}, \ \vec{w}, \ { \sf sign}( \vec{v}+(1|0)), 
	 \ { \sf sign}( \vec{w}+(1|0)), \ { \sf sign}( \vec{v}-(1|0)), \ { \sf sign}( \vec{w}-(1|0)) $$
	 are located on the same line \ $L$, and they have to be  elements of \ $Set1$ . \  
	 The point \ $ \vec{v}$ \ could be \ $ P_{hor}$ .  Then  \ 
	 $ S =   { \sf sign}( \vec{v}+(1|0))$ \ and \ $ T =  { \sf sign}( \vec{v}-(1|0)) $ .  If we imagine that
	 $ \vec{w}$ is located away from \ $ \vec{v}$ \ in one direction on \ $L$ \ or the other, either \ 
	 ${ \sf sign}( \vec{w}+(1|0))$ \ is not in \ $Set1$ \ or \ ${ \sf sign}( \vec{w}-(1|0))$ \ is not 
	 in \ $Set1$. Hence it is only possible that \ $ \vec{v} = \vec{w} = P_{hor}$ , \ which contradicts our 
	 assumption.  Thus we have proved   Case ({ \bf A1}).  \\  \\
	 Finally  follows the last case ({ \bf A2}), because if we  can prove this, the other cases 
	 ({\bf B1}) and \  ({\bf B2}) can be shown in the same manner, and no other ideas are  needed. \
	 Hence we only prove \\ { Case ({ \bf A2}). }  \
	 Let  \quad $   v_1 < w_1 $ \quad and \quad  $ \{v_1, w_1\} \cap [-1,+1] \neq  \emptyset $ ,
	 for instance let \ $ -1 \leq w_1 \leq 1 $ .        
	  The two different points \ $\vec{v}, \vec{w}$ \ define an  unique straight line \
	  $ L: {\tt y} = m_L \cdot {\tt x} + b_L \ , \ m_L,b_L \in  \mathbbm{R}$.                                                                          By  \ lemma   \ref{lemma F}  \  
	  	  all seven points \   $ (0|1),  \vec{v}, \ \vec{w}, \ { \sf sign}( \vec{v}+(1|0)), 
	 \ { \sf sign}( \vec{w}+(1|0)), \ { \sf sign}( \vec{v}-(1|0)), \\ { \sf sign}( \vec{w}-(1|0)) $ \ 
	 are located on \ $L$, \ hence \ $b_L = 1 $ .  
	 \begin{lemma}  \sf
	 We have that \ \ $ -1 < m_L < 1 $  .
	 \end{lemma}
	 \begin{proof}
	 If \ \ $ m_L < -1 $ \ \ or \ \ $ m_L > +1 $ , \ \ it would contradict the convexity of the unit ball \ $ {\bf B} $ .      (Note the four unit vectors \ $ (-1|0), \vec{w}, \vec{v}, (1|0) ) $. \ \
	 Now we assume  \  $m_L \in \{ -1,1 \} $ , \ for instance \  $m_L = 1 $ . \   
	 Then  \ $ \vec{w}, \vec{v} $ \ are on \ \ $ L: {\tt y} = {\tt x} + 1$  \ , hence the three vectors \ 
	 $ (0|0),  \vec{v}+(1|0), \vec{w}+(1|0) $ \  would be collinear on the line \ $  {\tt y} = {\tt x} $. \
	 Because of \  $  \left\|\vec{v}+(1|0)\right\| =  h_+(t_1) =  h_+(t_2) = \left\|\vec{w}+(1|0)\right\| $ \ 
	 it follows that \  $  \vec{v}+(1|0) = \vec{w}+(1|0) $ , \  hence \ $  \vec{v} = \vec{w}$ , \ and we get a                       contradiction.
	 \end{proof}
	  Now we use proposition \ref{proposition B} .  \ Abbreviate 
	  \ $ \overline{T}_w := { \sf sign}( \vec{w}+(1|0)), \ \overline{S}_w := { \sf sign}( \vec{w}-(1|0)), \ 
	   \overline{T}_v := { \sf sign}( \vec{v}+(1|0)), $ \ and \ $ \overline{S}_v := { \sf sign}( \vec{v}-(1|0)) $.  
	  Now consider the eight unit vectors \ $  \overline{T}_w, \overline{S}_w,  \overline{T}_v, \overline{S}_v,$
	  \ and their negatives\ $ -\overline{T}_w, -\overline{S}_w,  -\overline{T}_v, -\overline{S}_v $ .   
  Four at a time create a parallelogram, namely \ $\overline{T}_w, \overline{S}_w, -\overline{T}_w, -\overline{S}_w$,     \ and \ $ \overline{T}_v, \overline{S}_v, -\overline{T}_v, -\overline{S}_v$, respectively.  
    By  proposition \ref{proposition B} ,    
    the line that connects $ \overline{T}_v$ \ and \ $  -\overline{S}_v$ \ and 
    also the line that connects  $ \overline{T}_w $ \ and \ $ -\overline{S}_w$ \  meet the horizontal {\tt x}-axis \
    in the point \ $(1|0)$ ,  and, because of \ $ -1 < m_L < 1 $ ,   
     by proposition \ref{proposition B} 
     the two points $ \overline{T}_v , \overline{S}_v$ \ and  \  $ \overline{T}_w  , 
     \overline{S}_w$ , \ respectively, both are located  above 
     the horizontal axis, hence the intersection point \ $(1|0)$ \ is between \ $ \overline{T}_v \ , 
    \  -\overline{S}_v$ \ \ and \ \  $ \overline{T}_w  \ , \  -\overline{S}_w$ , respectively.  
     \  Now let us consider the line \ $J$ \ that connects the unit vectors \ 
     $  \overline{T}_w$ \ and \ $-\overline{S}_v$ .
      Because  of the convexity of  $ {\bf B}$, \ $J$ \  must be a subset of \ $ {\bf B}$. 
       Because  of  our  assumption \ $ v_1 < w_1 $  \  , \	 
   on the line $L$ \  the most left one of the six different unit vectors  \ \ 
    $ \overline{S}_v, \overline{S}_w, \vec{v}, \vec{w}, \overline{T}_v $ \ and \ $ \overline{T}_w$ \ \ is 
    $ \overline{S}_v$ , and the most right one of the six is \  \ $ \overline{T}_w$ .     
    Correspondingly,  of the  six  unit vectors 
     $ -\overline{S}_v, -\overline{S}_w, -\vec{v}, -\vec{w}, -\overline{T}_v, -\overline{T}_w$ , \ \   
    the most  left one of the six is \  $ -\overline{T}_w$ ,  and the most right one is \ $ -\overline{S}_v$ .    
   We take on the lines $L$ \ and \ $-L$, respectively, the points that are the most right ones, namely
    \ $  \overline{T}_w$ \ and \ $-\overline{S}_v$ .  \ The line \ $J$ \  that  connects both points 
     \ $  \overline{T}_w$ \ and \ $-\overline{S}_v$ \ crosses the {\tt x}-axis in \ $\lambda > 1$ . \
     Because of  $J \subset {\bf B}$,
      for the norm of $(\lambda|0)$ \ holds that \ $ \left\|(\lambda|0)\right\| \leq 1$ . \ 
     Because of \  \ \ $\lambda > 1 , \ (1|0) $ \ would be in the interior  of ${\bf B}$. 
       That contradicts \
     $ \left\|(1|0)\right\| = 1$ , \ and finally we have found a contradiction also for the last
      Case ({ \bf A2}) .  \\  That means that all cases ({ \bf A1}), ({ \bf A2}), ({ \bf B1}), ({ \bf B2}),
      ({ \bf C1}), ({ \bf C2}), and ({ \bf C3}) have been discussed,  and  with the assumption \  '$ t_1 < t_2 $' , 
      that means \ '$ \vec{v} \neq \vec{w}$' ,  \  
      we always  found a contradiction,  thus \ $ \vec{v} = \vec{w} $  \ and \ $ t_1 = t_2 $ \ 
      remains as  the only possibility.  Hence the map $ \Theta $ is injective, hence bijective,  
       and  finally  \ { \bf Theorem  \ref{main Theorem}  } \ has been proved \:! 
   \end{proof}
    \section{On Concave Corners and Some Open Problems}  
    In proposition  \ref{proposition1} \ we claimed and proved that in a  real seminormed vector space
     \ $( X , \|..\|) $ \ the 	 triple \ \
	      $( X , \|..\| , <.. \: | \: .. >_\spadesuit ) $ \  \  fulfils the  { \sf CSB } inequality, 
	      hence the  \ 'Thy angle' \  $ \angle_{Thy} (\vec{x}, \vec{y}) $ \  is defined for all  \
	       $ \vec{x}, \vec{y}  \neq \vec{0} $ .     
    Furthermore, in \ { \bf Theorem  \ref{main Theorem}  } \   we proved that if \ $(  X , \|..\| )$ \ even is
    a normed space,  the axiom   ({\tt An} 11) \ is fulfiled. For proving  both facts we always use the 
    convexity of the unit ball \  $ {\bf B}$  \ of  $(  X , \|..\| )$ . \  If you note that the 
    convexity of the unit ball is equivalent to the triangle inequality, 
    we come to the  natural question 
    whether there exists  a  \ {\bf hw \ space}  \ $ ( X , \|..\| \   ) $ , \  such that \ $\|..\|$  
     fulfils (1), the absolute homogenity, and (2), the positive definiteness, but not (3), the triangle inequality, 
    and \ $( X , \|..\| , <.. \: | \: .. >_\spadesuit ) $ \  satisfies the \  { \sf CSB } inequality,
     or even the axiom  ({\tt An} 11). 
    Natural candidates are the spaces \ $ ( \mathbbm{R}^{2} ,  \|..\|_{{ \tt Polygon},r}) $ , \
    which will  be  defined at once.    
    But in the case of   \ $ r < 1 $ \ we have no  success, as we will see. 
     For \ $r < 1$ , \  $ ( \mathbbm{R}^{2} ,  \|..\|_{{ \tt Polygon},r}) $ 
    has something that we now  call a ' concave corner'.    
    \begin{definition}               \rm
    \quad Let  the pair \ $ ( X , \|..\| \   ) $  \ \rm  be a  \ {\bf hw \ space}, \
                              let  \  $ \widehat{y} \in X $  \ , \ $ \| \widehat{y}\| > 0 $ . \   \\   
    $ \widehat{y} $ \  is  called  a  { \it concave   corner } \ \ $ \Longleftrightarrow $ \ \                     
    there is  an $   \overline{x} \in X $ , $ \| \overline{x}\| > 0 $ ,  \ and there are three real numbers  \
    \ $ \varepsilon, \ m_-, \ m_+ $ , \ with \ $   m_- < m_+ $ \ and  \  $ \varepsilon > 0 $ , \ such that 
    for all \ $ \delta \in [ 0, \varepsilon ] $ \ we have that \   
    $ \|  \delta \cdot { \sf sign}(\overline{x})+ ( 1 + \delta \cdot m_+ ) \cdot  {\sf sign}(\widehat{y}) \| =  1 
     =  \|  -\delta \cdot {\sf sign}( \overline{x}) + ( 1 - \delta \cdot m_- ) \cdot {\sf sign}(\widehat{y}) \| $ .
   \end{definition}
    We get a set of homogeneous weights on  \ $  \mathbbm{R}^{2}$ \ if we define  for every  \ \ $ r > 0 $ \
   a homogeneous weight  \ \ 
    $ \|..\|_{{ \tt Polygon},r} :   \mathbbm{R}^{2} \longrightarrow \mathbbm{R}^{+} \cup \{0\} $ , if 
   we fix the unit sphere  $ {\bf S}$ \  of  $ ( \mathbbm{R}^{2}, \|..\|_{{ \tt Polygon},r} )  $   with the polygon        through the six
   points \ $ \{ (0|r), (1|1), (1|-1), (0|-r),(-1|-1),(-1|1) \} $ \ and returning to $(0|r)$ ,
   and then extending \  $ \|..\|_{{ \tt Polygon},r}$ \ by  homogenity. \\  
     ( See  Picture 8) .   \\  \\    \\  
         \quad \\  
      \setlength{\parindent}{0mm}
    \begin{figure}[ht]
    \centering
     \setlength{\unitlength}{1cm}
    \begin{picture}(6,4)     \thinlines  
      \put(6.15,-0,1){$   {\tt x}   $}     \put(-0.3,5.8){$  {\tt y} $}    
       \put(-6.0,0){\vector(1,0){12.0}}     \put(0,-4){\vector(0,1){10}}  
               \put(-4.0,-4.0){\line(0,1){8.0}}   \put(4.0,4.0){\line(0,-1){8.0}} 
          \put(-4.0,4.0){\line(2,-1){4.0}}    \put(0.0,2.0){\line(2,1){4.0}}     \put(0.0,-2.0){ $ -\frac{1}{2} $ } 
    \put(0.0,1.7){ $ \frac{1}{2} $ }    \put(-0.5,1.7){ $ \widehat{y} $ }    \put(3.55,-0.35){ $ \overline{x} $ }   
            \put(-4.0,-4.0){\line(2,1){4.0}}   \put(0.0,-2.0){\line(2,-1){4.0}} 
                \put(0.0,4.0){ $ 1 $}      
                 \put(-0.27,4.0){ $ \times $}     \put(3.9,-0.36){ $ 1 $}  \put(-4.67,-0.36){ $ -1 $} 
    \put(0.5,5.2){ The unit sphere of  \ $ ( \mathbbm{R}^{2} ,  \|..\|_{{ \tt Polygon},\frac{1}{2}}) $ ,      
                   with the concave corner  }  
  \put(0.5,4.4){   $ \widehat{y} = (0|\frac{1}{2})$ . \quad 
                      We have \quad $ \overline{x} = (1|0) , \ \varepsilon = 1 , \ m_- = -1 < m_+ = +1 $ . }      
   \end{picture}
    \end{figure}	    
    \quad  \\  \\  \\  \\   \\  \\ \\  \begin{center}  Picture 8    \end{center} { \ }      
    \begin{lemma}     \sf
     For all \  $ 0 < r < 1 $ , \ the space \ \ $ ( \mathbbm{R}^{2} ,  \|..\|_{{ \tt Polygon},r}) $ \  has
     a  concave corner at \ $ \widehat{y} := (0|r) $, \ with  $ \overline{x} := (1|0) \ , \ \varepsilon := 1 \ 
             , \ m_- := 1 - \frac{1}{r} < 0 < m_+ :=  \frac{1}{r} - 1 $ .  
     \end{lemma} 
    \begin{proposition}  \qquad   \sf 
    \quad Let  the pair \ $ ( X , \|..\| \   ) $  \ \rm  be a  \ {\bf hw \ space}, \
          let  \  $ \widehat{y} \in X $  \ \ be a  concave   corner.   \ \   
    Then  \ $ ( X , \|..\|, <.. \: | \: .. >_\spadesuit ) $  \  does  not fulfil the   {\sf CSB } inequality.
     \end{proposition} 
   \begin{proof}  \quad   We use the above elements \ \ $ \widehat{y}, \ \overline{x} $ \ \ and  then  for all
     \ \ $ \delta \in [ 0, \varepsilon ] $ \ \ we can compute \\ $ P_\spadesuit(\delta) \ := \  
    <  \delta \cdot { \sf sign}(\overline{x})+ ( 1 + \delta \cdot m_+ ) \cdot  {\sf sign}(\widehat{y}) \ \: | 
    \:  -\delta \cdot {\sf sign}( \overline{x}) + ( 1 - \delta \cdot m_- ) \cdot {\sf sign}(\widehat{y}) >_\spadesuit \\
     = \ \    \frac{1}{4} \cdot   
    \left[  \ \left\| \ [ 2 + \delta \cdot ( m_+ - m_- )] \cdot {\sf sign}(\widehat{y}) \ \right\|^{2} \  - 
    \ \left\| \   2 \cdot    \delta \cdot {\sf sign}( \overline{x}) + \delta \cdot ( m_+ + m_- ) \cdot 
                                  {\sf sign}(\widehat{y}) \ \right\|^{2} \  \right]   \\
     = \ \    \frac{1}{4} \cdot   
    \left[  \  [ 2 + \delta \cdot ( m_+ - m_- )]^{2} \cdot  \left\| \ {\sf sign}(\widehat{y}) \ \right\|^{2} \  - 
    \ \delta^{2} \cdot \left\| \   2 \cdot {\sf sign}( \overline{x}) +  ( m_+ + m_- ) \cdot 
                                  {\sf sign}(\widehat{y}) \ \right\|^{2} \  \right]   \\
     = \ \   1 + \delta \cdot ( m_+ - m_- ) + \frac{1}{4} \cdot \delta^{2} \cdot \left[ \ ( m_+ - m_- )^{2} \ -  \        \left\| \ 2 \cdot {\sf sign}( \overline{x}) +  ( m_+ + m_- ) \cdot {\sf sign}(\widehat{y}) \ \right\|^{2} \ \right]        \\   = \ \   1 + \delta \cdot ( m_+ - m_- ) + \frac{1}{4} \cdot \delta^{2} \cdot K  $ \ ,  \\                
   with the real constant \ \ \ $  K :=  \ ( m_+ - m_- )^{2} \ -  \ \left\| \ 2 \cdot {\sf sign}( \overline{x}) + 
                                      ( m_+ + m_- ) \cdot {\sf sign}(\widehat{y}) \ \right\|^{2} $ . \\
   This calculation holds for all \ $ \delta \in   [ 0, \varepsilon ] $ . \  Hence, 
    for a positive \ $\widehat{\delta}$ that is almost 0 , because of \ $  m_+ - m_- > 0 $ \ we have  \ \ 
   $ P_\spadesuit(\widehat{\delta})  > 1 $ , hence, for the unit vectors \\  
    $  \widehat{\delta} \cdot { \sf sign}(\overline{x})+ ( 1 + \widehat{\delta} \cdot m_+ ) \cdot  
    {\sf sign}(\widehat{y}) $  \quad \  and  \quad \ $ -\widehat{\delta} \cdot {\sf sign}( \overline{x}) + 
    ( 1 - \widehat{\delta} \cdot m_- ) \cdot {\sf sign}(\widehat{y}) $ , \\ 
    the   { \sf CSB } inequality is not satisfied.  
   \end{proof}
   \begin{corollary}   \sf
   For all \  $ 0 < r < 1 $ , \ the space \ \ 
   $ ( \mathbbm{R}^{2} ,  \|..\|_{{ \tt Polygon},r} , <.. \: | \:..>_\spadesuit ) $ \ \ does not 
   fulfil the \  { \sf CSB } inequality. \ Hence, there are vectors \ $ \vec{x} \neq \vec{0} \neq \vec{y}$ \ 
   such that  the  \ 'Thy angle' \  $ \angle_{Thy} (\vec{x}, \vec{y}) $ \  is not defined. 
   Hence, for  \  $ 0 < r < 1 $ , \ the space \ \ 
   $ ( \mathbbm{R}^{2} ,  \|..\|_{{ \tt Polygon},r} ,  \angle_{Thy} ) $ \ \ is  not an 'angle space' .  
   \end{corollary} 
    As we mentioned in the beginning of the section, we formulate two open problems. \ We always deal  with  
     pseudonormed vector spaces   \ $ ( X , \|..\| \   ) $  \ ( that are \ {\bf hw \ spaces} \ which are 
   positive definite, that means 
     \ $ \|\vec{x}\| = 0 \ \ { \rm only \ for } \ \ \vec{x} = \vec{0} $ ) \  with a non-convex unit ball.  \\
    \underline{Problem 1:} \quad Exists  a  pseudonormed vector space \ $ ( X , \|..\| \   ) $ , \
            such that \ $\|..\|$  does 
     not fulfil  the triangle inequality,  but   \  \ $ ( X , \|..\|, <.. \: | \: .. >_\spadesuit ) $ 
      \ satisfies the \  { \sf CSB } inequality ? \\ 
    \underline{Problem 2:} \quad If \ {Problem 1} is true, \ exists  a  pseudonormed vector space \ $
                     ( X , \|..\| \   ) $ , \  such that \ $\|..\|$  
     does not fulfil  the triangle inequality, but \  \ $ ( X , \|..\|, <.. \: | \: .. >_\spadesuit ) $  \ 
      \ satisfies the \  { \sf CSB } inequality     and the axiom ({\tt An} 11) is fulfiled ? \\ 
    Good candidates for both questions are the \ { \it Hölder weights } $ \|..\|_p $ \ on   $\mathbbm{R}^{2}$ , \ \ 
   $ \|(x,y)\|_p \ = \ \sqrt[p]{ |x|^{p} + |y|^{p} }$ , \ \ with \ $ 0 < p < 1 $ , \  but  \ $p$ \ almost 1 .     
      \section{A Generalization of the  Thy Angle}  
       Recall in  the general  definitions the  expression  of a  convex set,  and   for a  subset $ A $ 
       of \ \ $ ( X ,  \|..\| )$  \  we defined   the convex hull \ { \it conv($A$) } ,  
  $$  { \it conv(A) } := \bigcup \ \{ \sum_{i=1}^{n} t_i \cdot \vec{x_i}  \ | \ n \in \mathbbm{N} , \ t_i \in [0,1] \ 
 \ {\rm and } \ \ \vec{x_i} \in A  \ \ { \rm for} \ i = 1, ... , n \ , \ {\rm and } \ \sum_{i=1}^{n} t_i = 1 \} \ , $$          which is the smallest convex set that contains $A$.  \\
    For a  	{\bf hw space} \ $( X , \|..\| )$ \ with the unit ball \ $ {\bf B} := \{ \: \vec{x} \in  X \ |  
	    \ \|\vec{x}\| \leq 1 \: \} $ , we defined a seminorm 	 $ \|..\|_{| \it conv({\bf B})}$ , for all 
	     \ $ \vec{x} \in X $, \ \ let  \ \
     $ \|\vec{x}\|_{| \it conv({\bf B})} := 
              \ \inf \: \{ r > 0 \: | \: \frac{1}{r} \cdot \vec{x} \in { \it conv({\bf B})} \} $.  \\
   \ We have  for all \ $ \vec{x} \in X $ \ that  \ $ \|\vec{x}\|_{| \it conv({\bf B})} \ \leq  \  \|\vec{x}\| $ . \  
   Note that \ \  $( X , \|..\| )  \stackrel{id}{\rightarrow} ( X , \|..\|_{| \it conv({\bf B})})$ \ \
    is a continuous map. Further    
     note that for a 	{\bf hw space} \ $( X , \|..\| )$ \ with the unit ball \ $ {\bf B} $ , \ 	the pair \ 
	 $( X , \|..\|_{| \it conv({\bf B})})$ \ is a seminormed vector space,  because  the triangle inequality
	  is satisfied.  \  
	   Then we called \ $ \|..\|$, \ or  the pair \ ( $ X , \|..\|$ ), \ respectively, \   { \it normable} \ \
   if and only if  the pair \  $( X , \|..\|_{| \it conv({\bf B})})$ \ is a normed vector space. \\
    In a  normable space \  ($ X , \|..\|$) \ we have  the zero-set \ $ { \sf Z} = \{ \vec{0} \} $ .  \ 
    Then there is an equivalence  \ \   $ \vec{x} \neq \vec{0}  \ \Longleftrightarrow \ 
     \|\vec{x}\|  > 0  \ \Longleftrightarrow \   \|\vec{x}\|_{| \it conv({\bf B})}  > 0 $ , \ for all \ 
    \ $ \vec{x} \in X $ . 
    
    \begin{definition}     \rm   
        \quad Let  the pair \ $ ( X , \|..\| \   ) $  \ \rm  be a  \ normable {\bf hw space}.
       We  define   a continuous product  \quad 
         $  < .. \: | \:  .. >_{\spadesuit \ | \it conv}  \ : \  X^{2} \longrightarrow   \mathbbm{R} $ .  \\
    If  \ \ \ $ \vec{x} = \vec{0} $ \ or \ $ \vec{y} = \vec{0} $  
      \ \ we  set \ \  $  < \vec{x} \: | \:  \vec{y} >_{\spadesuit \ | \it conv} \ := \ 0 $ ,  \ 
       and   in the  case of  \ \  $ \vec{x} \neq \vec{0} \neq \vec{y}  $  
         \quad we  define  \ \quad  $  < \vec{x} \: | \:  \vec{y} >_{\spadesuit \ | \it conv} \ \  :=  \ \  $         \begin{equation*}   
     \frac{1}{4} \cdot  \|\vec{x}\| \cdot \|\vec{y}\| \cdot 
         \left[  \ \left\| \frac{\vec{x}}{\|\vec{x}\|_{| \it conv({\bf B})}}   + 
          \frac{\vec{y}}{\|\vec{y}\|_{| \it conv({\bf B})}} \right\|_{| \it conv({\bf B})}^{2} \  - \
         \left\| \frac{\vec{x}}{\|\vec{x}\|_{| \it conv({\bf B})}}   -   
         \frac{\vec{y}}{\|\vec{y}\|_{| \it conv({\bf B})}} \right\|_{| \it conv({\bf B})}^{2} \  \right] .
   \end{equation*}   
     \end{definition}  
  \begin{proposition}     \sf   \quad
     Let \ \ $( X , \|..\| )$ \ \ be a normable \ {\bf hw space} \ with the unit ball \ $ {\bf B}$ . \\
   (1) \ \ The    product  \ \ 
     $  < .. | .. >_{\spadesuit \ | \it conv}  \ : \  X^{2} \longrightarrow   \mathbbm{R} $    
  \ \  fulfils  the properties \ \   $\overline{(1)}$  \  ("homogenity"),  \ $\overline{(2)}$  \ ("symmetry"), \               $\overline{(3)}$ \ ("positive semidefiniteness") , \ and  \ $\overline{(4)}$ \ ("definiteness") . \ \
    Hence,  \ \ $ ( X ,  < .. | .. >_{\spadesuit \ | \it conv} )$ \ \ is a homogeneous product vector space.  \\  
  (2) \ \ We have \ \ $ \|\vec{x}\| =  \sqrt{ <\vec{x} \: | \: \vec{x}>_{\spadesuit \ | \it conv} } $
                               \ \ for all $ \vec{x} \in X $ .                        \\ 
   (3) \ \ The triple \ \ $( X , \|..\|,  < .. | .. >_{\spadesuit \ | \it conv}   )$ \ \  fulfils the 
    { \sf CSB } inequality, that  means  for  all \ \ $ \vec{x} , \vec{y} \in X $ \ \ it holds that \ \ \
     $  |< \vec{x} \: | \:  \vec{y} >_{\spadesuit \ | \it conv}| \ \leq \ \|\vec{x}\|  \cdot    \|\vec{y}\| $ . \\
   (4) \ \ In the case of a normed space \ \  $( X , \|..\| ) $  \ \  we get that  \ \ 
   $ \|..\| =  \|..\|_{| \it conv({\bf B})}  $ , \ \  and \ \ $ < .. | .. >_{\spadesuit \ | \it conv}  \  
    = \  < .. | .. >_{\spadesuit} $ .      
  \end{proposition} 
    \begin{proof}  \rm   \quad 
  (1) \quad  See the comments after the definition of the product \ $   < .. | .. >_{\spadesuit} $ \ on page 5 .    \\
   (2) \quad  Trivial.   \\ 
   (3) \quad  Easy, because \ $ \|..\|_{| \it conv({\bf B})} $ \ fulfils the triangle inequality.  \\
   (4)  \quad  Trivial, \  because  \ {\bf B} \ is convex, \ hence \ {\bf B} = {\it conv}({\bf B})  .  
  \end{proof}  
  	\begin{definition}    \rm
	  For all   normable {\bf hw spaces}  ( $ X , \|..\|$ ), \   
	   for all  \ $ \vec{x}, \vec{y} \in X\backslash { \sf Z}  $ , \ ( that means \
	    $ \vec{x} \neq  \vec{0} $ \ and \ $ \vec{y}  \neq  \vec{0} $ ) \ \ with \  \
	    $ | < \vec{x} \: | \:  \vec{y} >_{\spadesuit \ | \it conv} | \leq \|\vec{x}\| \cdot \|\vec{y}\| $ , \
	  we define the   \\  { \it generalized  Thy angle }  \quad  \ 
	      $ \angle_{ { \bf \overline{Thy}} } (\vec{x}, \vec{y}) \ := \
  \arccos{\frac{ < \vec{x} \: | \:  \vec{y} >_{\spadesuit \ | \it conv}  }{ \|\vec{x}\| \cdot \|\vec{y}\| }} \ = $ 
       \begin{equation*} 
           \arccos{ \left( \frac{1}{4} \cdot   
        \left[  \ \left\| \frac{\vec{x}}{\|\vec{x}\|_{| \it conv({\bf B})}}   + 
          \frac{\vec{y}}{\|\vec{y}\|_{| \it conv({\bf B})}} \right\|_{| \it conv({\bf B})}^{2} \  - \ \
         \left\| \frac{\vec{x}}{\|\vec{x}\|_{| \it conv({\bf B})}}   -   
     \frac{\vec{y}}{\|\vec{y}\|_{| \it conv({\bf B})}} \right\|_{| \it conv({\bf B})}^{2} \  \right]   \right)} \ . 
    \end{equation*}    
		\end{definition} 
		\begin{proposition}  \qquad   \sf    \label{proposition34} 
		Let the pair \ ( $ X , \|..\|$ ) \ be a   normable {\bf hw space}.   \\ 
	  (a)    \ \  \  The 	 triple \ \
	      $( X , \|..\| , <.. \: | \: .. >_{\spadesuit \ | \it conv} ) $ \  \  fulfils the  { \sf CSB } inequality, 
	   hence the  generalized angle \ \  $ \angle_{ { \bf \overline{Thy}} } (\vec{x}, \vec{y}) $ \ \
	        is defined for all  \   $ \vec{x} \neq \vec{0} \neq \vec{y} $  .                        \\ 
	  (b)  \ \ \ The 	 triple \ 
	  \  $(X , \|..\|,\angle_{ { \bf \overline{Thy}} })$ \
	  fulfils  all the  demands   ({\tt An} 1), ({\tt An} 2), ({\tt An} 3), ({\tt An} 4), ({\tt An} 5), ({\tt An} 6),              ({\tt An} 7) . \ Hence  \ \
	   ( $ X , \|..\|,\angle_{ { \bf \overline{Thy}} } $ ) \ is  an angle space.  \\
	   (c)  \ \ \ If even the pair \ ( $ X , \|..\|$ ) \ is a  normed vector space,  then 
	   we have \    for all \  $ \vec{x},\vec{y} \neq \vec{0} $ \  that \ \   
	  $ \angle_{ { \bf \overline{Thy}} } (\vec{x}, \vec{y})  \  =  \ \angle_{Thy} (\vec{x}, \vec{y})  $ .  \\
	   (d)   \ \  \sf  If even \  $ (  X , <.. \: | \: .. >_{IP} )$  \ \ is an  {\bf IP space}, then 	  
	   we have   \   for all \  $ \vec{x},\vec{y} \neq \vec{0} $ \  that \ \   
	  $ \angle_{ { \bf \overline{Thy}} } (\vec{x}, \vec{y})  \  =  \ \angle_{Euclid} (\vec{x}, \vec{y})  $ .   \\
	  (e)\ \ \ If \ ( $ X , \|..\|$ ) is a normable {\bf hw space}, then the 	 triple \ 
  	\ $(X , \|..\|,\angle_{ { \bf \overline{Thy}} })$ \ generally does not  fulfil \
  	 ({\tt An} 8), ({\tt An} 9), ({\tt An} 10) .  \\
   (f)\ \ \ If \ ( $ X , \|..\|$ ) is a normable {\bf hw space}, then  \ 
	  \  $(X , \|..\|,\angle_{ { \bf \overline{Thy}} })$ \  satifies \ ({\tt An} 11) .   \\ 
	  \end{proposition} 
	   \begin{proof}  \rm   \quad 
   (a) \quad  The \  { \sf CSB } inequality  was shown in the previous proposition.   \\
   (b) \quad  Easy. \ \ See the proofs in \ proposition  \ref{proposition1}.    \\ 
   (c) \quad  Trivial,
         because \ \ {\bf B} = {\it conv}({\bf B}) , \ hence \ \ $ \|..\| = \|..\|_{| \it conv({\bf B})}$ .  \\
   (d)  \quad  Trivial, \  because in an  {\bf IP space} is \ \ $ \angle_{Euclid} =  \angle_{Thy}$ .  \\
   (e)  \quad  Use the examples of \ proposition  \ref{proposition1}.    \\
   (f)  \quad  The proof is not trivial,  but the same as in \ { \bf Theorem  \ref{main Theorem}  } .  
   \end{proof}       
	  \quad  \\  \\  	 
	  { \bf Acknowledgements } \\ We wish to thank   \ \ 
     Dr. Gencho Skordev  and  Prof. Dr. Eberhard  Oeljeklaus ,  \ which supported us by 
     patient listening and interested discussions.  
	 \quad   \\   \\

	\end{document}